\documentclass[11pt]{article}
\usepackage{latexsym}
\usepackage[all]{xy}
\usepackage{amsmath}
\usepackage{amssymb}
\usepackage{amsfonts}

\newtheorem{thm}{Theorem}[section]
\newtheorem{prop}[thm]{Proposition}
\newtheorem{lem}[thm]{Lemma}
\newtheorem{sublem}[thm]{Sub-lemma}
\newtheorem{df}[thm]{Definition}
\newtheorem{cor}[thm]{Corollary}

\newtheorem{rmk}[thm]{Remark}

\begin{document}

\title{\textbf{The homotopy theory of
dg-categories and derived Morita theory}}
\bigskip
\bigskip

\author{\bigskip\\
Bertrand To\"en \\
\small{Laboratoire Emile Picard UMR CNRS 5580} \\
\small{Universit\'{e} Paul Sabatier, Bat 1R2} \\
\small{Toulouse Cedex 9}
\small{France}}

\date{September 2006}

\maketitle

\begin{abstract}
The main purpose of this work is to study the homotopy theory
of dg-categories up to quasi-equivalences.
Our main result is a description
of the mapping spaces between two dg-categories $C$ and $D$ in terms
of the nerve of a certain category of $(C,D)$-bimodules. We also prove that the
homotopy category $Ho(dg-Cat)$ possesses
internal Hom's relative to the (derived) tensor product
of dg-categories. We use these
two results in order to prove a derived version of Morita theory,
describing the morphisms between
dg-categories of modules over two dg-categories $C$ and $D$ as
the dg-category of $(C,D)$-bi-modules. Finally,
we give three applications of our results. The first one expresses Hochschild
cohomology as endomorphisms of the identity functor, as well as
higher homotopy groups of the
\emph{classifying space of dg-categories} (i.e. the nerve of
the category of dg-categories and quasi-equivalences between them).
The second application is the existence of a good theory of localization for
dg-categories, defined in terms of a natural universal property. Our last application
states that the dg-category of (continuous) morphisms between the dg-categories of
quasi-coherent (resp. perfect) complexes on two schemes (resp. smooth and proper schemes) is 
quasi-equivalent to the
dg-category of quasi-coherent (resp. perfect) complexes on their product.
\end{abstract}

\medskip

\tableofcontents

\bigskip

\section{Introduction}

Let $A$ and $B$ be two associative algebras (over some field $k$), and
$A-Mod$ and $B-Mod$ be their categories of right modules. It is well known that
any functor $A-Mod \longrightarrow B-Mod$ which commutes with colimits is of the form
$$\begin{array}{ccc}
A-Mod & \longrightarrow & B-Mod \\
M & \mapsto & M \otimes_{A} P,
\end{array}$$
for some $A^{op}\otimes B$-module $P$. More generally, there exists a natural equivalence
of categories between $(A^{op}\otimes B)-Mod$ and the category of
all colimit preserving functors $A-Mod \longrightarrow B-Mod$. This is known
as Morita theory for rings.

Now, let $A$ and $B$ be two
associative dg-algebras (say over some field $k$), together with their triangulated
derived category of right (unbounded) dg-modules $D(A)$ and $D(B)$.
A natural way of constructing triangulated functors
from $D(A)$ to $D(B)$ is by choosing $P$ a left
$A^{op}\otimes B$-dg-module, and considering the derived functor
$$\begin{array}{ccc}
D(A) & \longrightarrow & D(B) \\
M & \mapsto & M\otimes_{A}^{\mathbb{L}}P.
\end{array}$$
However, it is well known that there exist
triangulated functors $D(A) \longrightarrow D(B)$ that do not
arise from a $A^{op}\otimes B$-dg-module (see e.g. \cite[2.5, 6.8]{ds}). The situation is even
worse, as the functor
$$D(A^{op}\otimes B) \longrightarrow Hom_{tr}(D(A),D(B))$$
is not expected to be reasonable in any sense
as the right hand side simply does not possess a natural triangulated structure.
Therefore,
triangulated categories do not
appear as the right object to consider if one is looking for
an extension of Morita theory to dg-algebras.
The main purpose of this work is to provide
a solution to this problem by replacing
the notion of triangulated categories by the notion of
dg-categories. \\

\begin{center} \textit{dg-Categories} \end{center}

A dg-category is a category which is enriched over the monoidal category of
complexes over some base ring $k$. It consists of a set of objects
together with complexes $C(x,y)$ for two any objects $x$ and $y$, and
composition morphisms $C(x,y)\otimes C(y,z) \longrightarrow C(x,z)$ (assumed to be
associative and unital). As linear categories can be
understood as \emph{rings with several objects}, dg-categories can be thought as
\emph{dg-algebras with several objects}, the precise statement being that
dg-algebras are exactly dg-categories having a unique object.

From a dg-category
$C$ one can form a genuine category $[C]$ by keeping the same set of objects
and defining the set of morphisms between $x$ and $y$ in $[C]$ to be
$H^{0}(C(x,y))$. In turns out that a lot of triangulated categories
appearing in geometric contexts are of the form $[C]$ for some
natural dg-category $C$ (this is for example the case for
the derived category of a reasonable abelian category, as well as for the derived
category of dg-modules over some dg-algebra).
The new feature of
dg-categories is the notion of \emph{quasi-equivalences}, a mixture between
quasi-isomorphisms and categorical equivalences and which turns out to be the
right notion of equivalences between
dg-categories. Precisely, a morphism $f : C \longrightarrow D$ between two
dg-categories is a quasi-equivalence if it satisfies the following two conditions
\begin{itemize}
\item For any objects $x$ and $y$ in $C$ the induced morphism
$C(x,y) \longrightarrow D(f(x),f(y))$ is a quasi-isomorphism.
\item The induced functor $[C] \longrightarrow [D]$
is an equivalence of categories.
\end{itemize}

In practice we are only interested in dg-categories up to quasi-equivalences, and the main
object of study is thus the localized category $Ho(dg-Cat)$ of dg-categories
with respect to quasi-equivalences, or
better its refined simplicial version $L(dg-Cat)$
of Dwyer and Kan (see \cite{dk2}). The main purpose of this paper
is to study the simplicial category $L(dg-Cat)$,
and to show that a derived version
of Morita theory can be extracted from it. The key tool for us will be the existence
of a model structure on the category of dg-categories (see \cite{tab}), which will
allow us to use standard constructions of homotopical algebra (mapping spaces,
homotopy limits and colimits \dots) in order
to describe $L(dg-Cat)$. \\

\begin{center} \textit{Statement of the results} \end{center}

Let $C$ and $D$ be two dg-categories, considered as objects in $L(dg-Cat)$.
A first invariant is the homotopy type of the simplicial set of morphism
$L(dg-Cat)(C,D)$, which is well known to be weakly equivalent to the mapping space
$Map(C,D)$ computed in the model category of dg-categories (see \cite{dk,dk2}).
From $C$ and $D$ one can
form the tensor product $C\otimes D^{op}$ (suitably derived if necessary),
as well as the category $(C\otimes D^{op})-Mod$ of $C\otimes D^{op}$-modules
(these are enriched functors from $C\otimes D^{op}$ to the category of complexes).
There exists an obvious notion of quasi-isomorphism between $C\otimes D^{op}$-modules,
and thus a homotopy category $Ho((C\otimes D^{op})-Mod)$. Finally, inside
$Ho((C\otimes D^{op})-Mod)$ is a certain full sub-category of
\emph{right quasi-representable objects}, consisting of modules $F$ such that
for any $x\in C$ the induced $D^{op}$-module $F(x,-)$ is quasi-isomorphic to a
$D^{op}$-module of the form $D(-,y)$ for some $y\in D$ (see \S 3 for details).
One can then consider the category $\mathcal{F}(C,D)$ consisting of all
right quasi-representable $C\otimes D^{op}$-modules and quasi-isomorphisms between them.
The main result of this work is the following.

\begin{thm}\label{ti1}{(See Thm. \ref{t1})}
There exists a natural weak equivalence of simplicial sets
$$Map(C,D)\simeq N(\mathcal{F}(C,D))$$
where $N(\mathcal{F}(C,D))$ is the nerve of the
category $\mathcal{F}(C,D)$.
\end{thm}

We would like to mention that this theorem does not simply follow from the
existence of the model structure on dg-categories. Indeed, this model structure
is not simplicially enriched (even in some weak sense, as
the model category of complexes is for example), and there is no obvious
manner to compute the mapping spaces $Map(C,D)$. \\

As an important corollary one gets the following result.

\begin{cor}\label{cti1}
\begin{enumerate}
\item There is a natural bijection between $[C,D]$, the set of morphisms
between $C$ and $D$ in $Ho(dg-Cat)$, and
the isomorphism classes of right quasi-representable objects in
$Ho((C\otimes D^{op})-Mod)$.
\item For two morphism $f,g : C\longrightarrow D$ there is a natural weak equivalence
$$\Omega_{f,g}Map(C,D) \simeq Map(\phi(f),\phi(g))$$
where $Map(\phi(f),\phi(g))$ is the mapping space between the  $C\otimes D^{op}$-modules corresponding to $f$ and $g$.
\end{enumerate}
\end{cor}

The tensor product of dg-categories, suitably derived, induces a
symmetric monoidal structure on $Ho(dg-Cat)$. Our second main result states
that this monoidal structure is closed.

\begin{thm}\label{ti2}{(See Thm. \ref{t2})}
The symmetric monoidal category $Ho(dg-Cat)$ is closed. More precisely,
for any three dg-categories $A$, $B$ and $C$, there exists
a dg-category $\mathbb{R}\underline{Hom}(B,C)$ and
functorial isomorphisms in $Ho(SSet)$
$$Map(A,\mathbb{R}\underline{Hom}(B,C))\simeq Map(A\otimes^{\mathbb{L}} B,C).$$
Furthermore, $\mathbb{R}\underline{Hom}(B,C)$ is naturally
isomorphic in $Ho(dg-Cat)$ to the dg-category of
cofibrant right quasi-representable $B\otimes C^{op}$-modules.
\end{thm}

Finally, Morita theory can be expressed in the following terms. Let
us use the notation
$\widehat{C}:=\mathbb{R}\underline{Hom}(C^{op},Int(C(k)))$, where
$Int(C(k))$ is the dg-category of cofibrant complexes. Note that by our theorem
\ref{ti2} $\widehat{C}$ is also quasi-equivalent to the dg-category of
cofibrant $C^{op}$-modules.

\begin{thm}\label{ti3}{(See Thm. \ref{t3} and Cor. \ref{ct3})}
There exists a natural isomorphism in $Ho(dg-Cat)$
$$\mathbb{R}\underline{Hom}_{c}(\widehat{C},\widehat{D})\simeq
\widehat{C^{op}\otimes^{\mathbb{L}} D},$$
where $\mathbb{R}\underline{Hom}_{c}(\widehat{C},\widehat{D})$ is the full sub-dg-category
of $\mathbb{R}\underline{Hom}(\widehat{C},\widehat{D})$ consisting of morphisms commuting with infinite
direct sums.
\end{thm}

As a corollary we obtain the following result.

\begin{cor}\label{cti3}
There is natural bijection between
$[\widehat{C},\widehat{D}]_{c}$, the sub-set of $[\widehat{C},\widehat{D}]$ consisting
of direct sums preserving morphisms, and
the isomorphism classes in $Ho((C\otimes^{\mathbb{L}} D^{op})-Mod)$.
\end{cor}

\bigskip

\begin{center} \textit{Three applications} \end{center}

\smallskip

We will give three applications of our general results. The first
one is a description of the homotopy groups
of the classifying space of dg-categories $|dg-Cat|$, defined
as the nerve of the category of quasi-equivalences between dg-categories.
For this, recall that the Hochschild cohomology of a dg-category
$C$ is defined by
$$\mathbb{HH}^{i}:=[C,C[i]]_{C\otimes^{\mathbb{L}}C^{op}-Mod},$$
where $C$ is the $C\otimes^{\mathbb{L}}C^{op}$-module
sending $(x,y)\in C\otimes C^{op}$ to $C(y,x)$.

\begin{cor}{(See Cor. \ref{chh'}, \ref{crpic})}
\begin{enumerate}
For any dg-category $C$ one has
\item
$$\mathbb{HH}^{*}(C)\simeq
H^{*}(\mathbb{R}\underline{Hom}(C,C)(Id,Id)).$$
\item
$$\pi_{i}(|dg-Cat|,C)\simeq \mathbb{HH}^{2-i}(C) \qquad\forall i>2.$$
\item $$\pi_{2}((|dg-Cat|,C)\simeq Aut_{Ho(C\otimes C^{op}-Mod)}(C)\simeq
\mathbb{HH}^{0}(C)^{*}$$
\item $$\pi_{1}(|dg-Cat|,\widehat{BA})\simeq RPic(A),$$
where $A$ is a dg-algebra, $BA$ the dg-category
with a unique object and $A$ as its endomorphism, and
where $RPic(A)$ is the derived Picard group
of $A$ as defined for example in \cite{rz,ke2,yek}.
\end{enumerate}
\end{cor}

Our second application is the existence of localization for
dg-categories. For this, let $C$ be any
dg-category and $S$ be a set of morphisms in $[C]$.
For any dg-category $D$ we define
$Map_{S}(C,D)$ as the sub-simplicial set of
$Map(C,D)$ consisting of morphisms sending
$S$ to isomorphisms in $[D]$.

\begin{cor}{(See Cor. \ref{cloc})}
The $Ho(SSet_{\mathbb{U}})$-enriched functor
$$Map_{S}(C,-) : Ho(dg-Cat_{\mathbb{U}}) \longrightarrow
Ho(SSet_{\mathbb{U}})$$
is co-represented by an object $L_{S}(C) \in Ho(dg-Cat_{\mathbb{U}})$.
\end{cor}

Our final application will provide a proof of the
following fact, which can be considered as a possible answer to a folklore question
to know whether or not all triangulated functors between derived categories
of varieties are induced by some object in the derived category of their product
(see e.g. \cite{o} where this is proved for
triangulated equivalences between derived categories of
smooth projective varieties).

\begin{cor}\label{cti3'}{(See Thm. \ref{tfour})}
Let $X$ and $Y$ be two quasi-compact and separated $k$-schemes, one of them
being flat over $Spec\, k$, and let $L_{qcoh}(X)$ and
$L_{qcoh}(Y)$ their dg-categories of
(fibrant) quasi-coherent complexes. Then, one has a natural
isomorphism in $Ho(dg-Cat)$
$$L_{qcoh}(X\times_{k} Y)\simeq \mathbb{R}\underline{Hom}_{c}(L_{qcoh}(X),L_{qcoh}(Y)).$$
In particular, there is a natural bijection between
$[L_{qcoh}(X),L_{qcoh}(Y)]_{c}$ and set of isomorphism classes of
objects in the category $D_{qcoh}(X\times Y)$.

If furthermore $X$ and $Y$ are smooth and proper over $Spec\, k$,
then one has
a natural
isomorphism in $Ho(dg-Cat)$
$$L_{parf}(X\times_{k} Y)\simeq
\mathbb{R}\underline{Hom}(L_{parf}(X),L_{parf}(Y)),$$
where $L_{parf}(X)$ (resp. $L_{parf}(Y)$) is the full sub-dg-category
of $L_{qcoh}(X)$ (resp. of $L_{qcoh}(Y)$) consisting of
perfect complexes.
\end{cor}

\bigskip

\begin{center} \textit{Related works} \end{center}

The fact that dg-categories provide natural and
interesting enhancement of derived categories
has been recognized for some times,
and in particular in \cite{bk}. They have been
used more recently in \cite{bll} in which a very special
case of our theorem \ref{tfour} is proved
for smooth projective varieties.
The present work
follows the same philosophy that dg-categories are
the \emph{true derived categories} (though I do not like
very much this expression).

Derived equivalences between (non-dg) algebras have been
heavily studied by J. Rickard (see e.g. \cite{ri1,ri2}),
and the results obtained have been commonly called
\emph{Morita theory for derived categories}.
The present work can be considered as
a continuation of this fundamental work, though our techniques
and our purposes are rather different. Indeed, in our mind
the word \emph{derived} appearing in our
title does not refer to generalizing Morita theory from module categories to
derived categories, but to generalizing
Morita theory from algebras to dg-algebras.

Morita theory for dg-algebras and ring spectra
has been approached recently using
model category techniques in \cite{ss}. The results
obtained this way state in particular that
two ring spectra have Quillen equivalent
model categories of modules if and only if
a certain bi-module exists. This approach, however, does
not say anything about \emph{higher homotopies}, in the
sense that it seems hard (or even impossible) to compare
the whole model category of bi-modules 
with the category of Quillen equivalences, already simply
because a model category of Quillen functors does
not seem to exist in any reasonable sense. This is another
incarnation of the principle that
model category theory does not work very well
as soon as categories of functors are involved, and
that some sort of higher categorical structures
are then often needed
(see e.g. \cite[\S 1]{to2}). 

A relation between the derived Picard group
and Hochschild cohomology is given in \cite{ke2}, and
is somehow close to our Corollary \ref{chh'}. An interpretation
of Hochschild cohomology as first order deformations of
dg-categories is also given in \cite{hagII}.

There has been many works on dg-categories (as well as its weakened, but
after all equivalent, notion
of $A_{\infty}$-categories) in which several universal constructions, such as
reasonable dg-categories of dg-functors or quotient and localization
of dg-categories, have been studied (see for example \cite{dr,ke,ly,ly2}). Of course,
when compared in a correct way, our constructions give back the same objects as the ones
considered in these papers, but I would like to point out that the two approaches
are different
and that our results can not be deduced from these previous works. 
Indeed, the universal properties of the constructions of
\cite{dr,ke,ly,ly2} are expressed in a somehow
un-satisfactory manner (at least for my personal taste)
as they are stated in terms of certain dg-categories of dg-functors
that are not themselves defined by some universal properties (except an obvious one
with respect to themselves !)\footnote{The situation is very comparable to the situation where one tries
to explain why categories of functors give the \emph{right notion}: expressing universal
properties using itself categories of functors is not helpful.}.
In some sense, the results proved in these papers are more
properties satisfied by certain constructions rather than existence theorems.
On the contrary our results truly are existence theorems and our dg-categories of dg-functors,
or our localized dg-categories,
are constructed as solution to a universal problem inside the category
$Ho(dg-Cat)$ (or rather inside the simplicial category $L(dg-Cat)$).
As far as I know, these
universal properties were not known to be satisfied by the constructions of
\cite{dr,ke,ly,ly2}.

The results of the present work can also be generalized in an obvious way to other
contexts, as for example simplicially enriched categories, or even
spectral categories. Indeed, the key tool that makes the proofs working is
the existence of a nice model category structure on enriched categories.
For simplicial categories this model structure is known to exist by a recent
work of J. Bergner, and our theorems \ref{t1} and \ref{t2} can be easily shown to
be true in this setting (essentially the same proofs work). Theorem
\ref{t3} also stays correct for simplicial categories except that one
needs to replace the notion of continuous morphisms by the more elaborated notion
of colimit preserving morphisms. More recently,
J. Tapia has done some progress
for proving the existence of a model category structure on $M$-enriched categories
for very general monoidal model categories $M$, including for example
spectral categories (i.e. categories enriched in symmetric spectra). I am convinced
that theorems \ref{t1} and \ref{t2}, as well
as the correct modification of theorem \ref{t3}, stay correct in this general setting.
As a consequence one would get a Morita theory for symmetric ring spectra.

Finally, I did not investigate at all the question of the behavior of
the equivalence of theorem \ref{t1} with respect to composition of morphisms. Of course,
on the level of bi-modules composition is given by the tensor product, but
the combinatorics of these compositions are not an easy question. This is related to the
question: \emph{What do dg-categories form ?} It is commonly expected that the
answer is \emph{an $E_{2}$-category}, whatever this means. The point of view of this work
is to avoid this difficulty by stating that
another possible answer is \emph{a simplicially enriched category} (precisely
the Dwyer-Kan localization $L(dg-Cat)$), which is a perfectly well understood structure.
Our theorem \ref{t2}, as well as its corollary \ref{cp5'} state that
the simplicial category $L(dg-Cat)$ is enriched over itself
in a rather strong sense. In fact, one can show that
$L(dg-Cat)$ is a \emph{symmetric monoidal simplicial category} in the
sense of Segal monoids explained in \cite{kt}, and I believe that
another equivalent way to talk about
$E_{2}$-categories is by considering $L(dg-Cat)$-enriched simplicial categories, again
in some Segal style of definitions (see for example \cite{to}). In other words,
I think the $E_{2}$-category of dg-categories should be completely determined by the
symmetric monoidal simplicial category $L(dg-Cat)$.

\bigskip

\textbf{Acknowledgments:} I am very grateful to M. Anel, C. Barwick, L. Katzarkov, T. Pantev,
M. Spitzweck, J. Tapia, M. Vaqui\'e and G. Vezzosi for their participation in
the small workshop on non-abelian Hodge theory which took place in Toulouse during
the spring 2004.
I also would like to thank warmly C. Simpson for his
participation to this workshop via
some tricky but enjoyable video-conference meeting. It has been during one of the
informal conversation of this workshop that the general ideas for a proof
of theorem \ref{tfour} have been found, and I think this
particular theorem should be attributed to the all of us.

I would like to thank B. Keller for several comments
on earlier versions of this work, and for pointing out to
me related references. I also thank the referee for
his careful reading.

\bigskip
\bigskip

\textbf{Conventions:}
All along this work universes will be denoted
by $\mathbb{U}\in \mathbb{V}\in \mathbb{W} \dots$.
We will always assume that they satisfy the
infinite axiom.

We use the notion of model categories
in the sense of \cite{ho}. The expression
\emph{equivalence} always refer to
weak equivalence in a model category.
For a model category $M$, we will denote by
$Map_{M}$ (or $Map$ if $M$ is clear) its mapping spaces
as defined in \cite{ho}. We will always consider
$Map_{M}(x,y)$ as an object in the homotopy category
$Ho(SSet)$. In the same way, the set of morphisms
in the homotopy category $Ho(M)$ will be denoted
by $[-,-]_{M}$, or by $[-,-]$ if $M$ is clear.
The natural $Ho(SSet)$-tensor structure
on $Ho(M)$ will be denoted
by $K\otimes^{\mathbb{L}}X$, for $K$ a simplicial set
and $X$ an object in $M$. In the same way, the
$Ho(SSet)$-cotensor structure will be denoted
by $X^{\mathbb{R}K}$. The homotopy fiber products will
be denoted by $x\times^{h}_{z}y$, and dually the
homotopy push-outs will be denoted by
$x\coprod^{\mathbb{L}}_{z}y$.

For all along this work, we fix an
associative, unital and commutative ring $k$. We denote
by $C(k)_{\mathbb{U}}$ the category of $\mathbb{U}$-small
(un-bounded) complexes of $k$-modules, for some
universe $\mathbb{U}$ with $k\in \mathbb{U}$.
The category $C(k)_{\mathbb{U}}$ is a
symmetric monoidal model category, where one uses
the projective model structures for which fibrations
are epimorphisms and equivalences are quasi-isomorphisms
(see e.g. \cite{ho}).
When the universe $\mathbb{U}$
is irrelevant we will simply write $C(k)$ for
$C(k)_{\mathbb{U}}$. The monoidal structure on $C(k)$ is the
usual tensor product of complexes over $k$,  and will be denoted by
$\otimes$. Its derived version will be denoted by $\otimes^{\mathbb{L}}$.

\section{The model structure}

Recall that a $\mathbb{U}$-small $dg$-category $C$ consists
of the following data.

\begin{itemize}
\item A $\mathbb{U}$-small set of objects
$Ob(C)$, also sometimes denoted by $C$ itself.
\item For any pair of objects
$(x,y)\in Ob(C)^{2}$ a complex
$C(x,y)\in C(k)$.

\item For any triple $(x,y,z)\in Ob(C)^{3}$
a composition morphism
$C(x,y)\otimes C(y,z) \longrightarrow C(x,z)$,
satisfying the usual associativity
condition.

\item For any object $x\in Ob(C)$, a
morphism $k \longrightarrow C(x,x)$, satisfying the
usual unit condition with respect to the
above composition.

\end{itemize}

For two dg-categories $C$ and $D$, a
morphism of dg-categories (or simply a dg-functor) 
$f : C \longrightarrow D$ consists
of the following data.

\begin{itemize}
\item A map of sets $f : Ob(C) \longrightarrow Ob(D)$.

\item For any pair of objects $(x,y)\in Ob(C)^{2}$,
a morphism in $C(k)$
$$f_{x,y} : C(x,y) \longrightarrow D(f(x),f(y))$$
satisfying the usual unit and associativity conditions.

\end{itemize}

The $\mathbb{U}$-small dg-categories and
dg-functors do form a category $dg-Cat_{\mathbb{U}}$.
When the universe $\mathbb{U}$ is irrelevant, we will simply write
$dg-Cat$ for $dg-Cat_{\mathbb{U}}$.
\\

We define a functor
$$[-] : dg-Cat_{\mathbb{U}} \longrightarrow Cat_{\mathbb{U}},$$
from $dg-Cat_{\mathbb{U}}$ to the category of $\mathbb{U}$-small
categories by the following construction.
For $C\in dg-Cat_{\mathbb{U}}$, the set of object of $[C]$
is simply the set of object of $C$. For two
object $x$ and $y$ in $[C]$, the set of morphisms
from $x$ to $y$ in $[C]$
is defined by
$$[C](x,y):=H^{0}(C(x,y)).$$
Composition of morphisms in $[C]$ is given by the natural
morphism
$$[C](x,y)\times [C](y,z) =
H^{0}(C(x,y)) \times H^{0}(C(x,y)) \longrightarrow
H^{0}(C(x,y)\otimes^{\mathbb{L}}C(y,z)) \longrightarrow
H^{0}(C(x,z))=[C](x,z).$$
The unit of an object $x$ in $[C]$ is simply given by
the point in $[k,C(x,x)]=H^{0}(C(x,x))$ image of the unit morphism
$k \longrightarrow C(x,x)$ in $M$.
This construction, provides a functor
$C \mapsto [C]$ from $dg-Cat_{\mathbb{U}}$ to the category
of $\mathbb{U}$-small categories. For a morphism $f : C \longrightarrow D$
in $dg-Cat$, we will denote by $[f] : [C] \longrightarrow [D]$ 
the corresponding morphism in $Cat$.

\begin{df}\label{d1}
Let $f : C \longrightarrow D$ be a morphism in
$dg-Cat$.
\begin{enumerate}
\item The morphism $f$ is \emph{quasi-fully faithful}
if for any two objects $x$ and $y$ in  $C$
the morphism $f_{x,y} : C(x,y) \longrightarrow D(f(x),f(y))$
is a quasi-isomorphism.
\item The morphism $f$ is \emph{quasi-essentially surjective}
if the induced functor $[f] : [C] \longrightarrow [D]$
is essentially surjective.
\item The morphism $f$ is a \emph{quasi-equivalence}
if it is quasi-fully faithful and quasi-essentially surjective.
\item The morphism $f$ is a \emph{fibration} if
it satisfies the following two conditions.
\begin{enumerate}
\item For any $x$ and $y$ in $C$ the morphism
$f_{x,y} : C(x,y) \longrightarrow D(f(x),f(y))$
is a fibration in $C(k)$ (i.e. is an epimorphism).
\item  For any $x\in C$, and any isomorphism
$v : [f](x) \rightarrow y'$ in $[D]$, there
exists an isomorphism $u : x \rightarrow y$ in $[C]$
such that $[f](u)=v$.
\end{enumerate}
\end{enumerate}
\end{df}

In \cite{tab} it is proved that
the above notions of fibrations and quasi-equivalences in
$dg-Cat$ form a model category structure. The
model category $dg-Cat_{\mathbb{U}}$ is furthermore
$\mathbb{U}$-cofibrantly generated in the sense
of \cite[Appendix]{hagI}. Moreover,
for $\mathbb{U}\in \mathbb{V}$, the set of
generators for the cofibrations and trivial cofibrations
can be chosen to be the same for $dg-Cat_{\mathbb{U}}$ and for
$dg-Cat_{\mathbb{V}}$. As a consequence we get that
the natural inclusion functor
$$Ho(dg-Cat_{\mathbb{U}}) \longrightarrow Ho(dg-Cat_{\mathbb{V}})$$
is fully faithful. This inclusion functor also induces
natural equivalences on mapping spaces
$$Map_{dg-Cat_{\mathbb{U}}}(C,D)\simeq
Map_{dg-Cat_{\mathbb{V}}}(C,D),$$
for two $\mathbb{U}$-small dg-categories $C$ and $D$.
As a consequence we see that we can change our universe
without any serious harm.

Note also that the functor
$$[-] : dg-Cat \longrightarrow Cat$$
induces a functor
$$Ho(dg-Cat) \longrightarrow Ho(Cat),$$
where $Ho(Cat)$ is the category of small categories and
isomorphism classes of functors between them. In other words,
any morphism $C \rightarrow D$ in $Ho(dg-Cat)$ induces
a functor $[C] \rightarrow [D]$ well defined up to a
non-unique isomorphism. This lack of uniqueness will not
be so much of a trouble as we will essentially be interested
in properties of functors which are invariant
by isomorphisms (e.g. being fully faithful, being an equivalence \dots). 

\begin{df}\label{dqess}
Let $f : C \longrightarrow D$ be a morphism of dg-categories.
The \emph{quasi-essential image of $f$} is the full sub-dg-category
of $D$ consisting of all objects $x\in D$ whose image in
$[D]$ lies in the essential image of the functor
$[f] : [C] \rightarrow [D]$.
\end{df}

The model category $dg-Cat$
also satisfies the following additional properties.

\begin{prop}\label{p0}
\begin{enumerate}
\item Any object $C\in dg-Cat$ is fibrant.

\item There exists a cofibrant replacement functor
$Q$ on $dg-Cat$, such that for any $C\in dg-Cat$ the natural
morphism $Q(C) \longrightarrow C$ induces the identity of the
sets of objects.

\item If $C$ is a cofibrant object in $dg-Cat$ and
$x$ and $y$ are two objects in $C$, then
$C(x,y)$ is a cofibrant object in $C(k)$.

\end{enumerate}
\end{prop}

\textit{Sketch of proof:} $(1)$ is clear by definition. $(2)$ simply follows from the fact that
one can choose the generating cofibrations $A\rightarrow B$ to induce the identity on the set
of objects (see  \cite{tab} for details). Finally, for $(3)$, one uses that any cofibrant object
can be written as a transfinite composition of push-outs along the generating
cofibrations. As the functor $C \mapsto C(x,y)$ commutes with filtered colimits, and that
a filtered colimit of cofibrations stays a cofibration, one sees that it is enough to
prove that the property $(3)$ is preserved by push-outs along a generating cofibration.
But this can be easily checked by an explicit description of such a push-out
(see \cite{tab} proof of Lem. 2.2. for more details). \hfill $\Box$ \\

To finish this paragraph, recall that a morphism
$x \rightarrow y$ in a model category $M$ is called a \emph{homotopy monomorphism}
if for any $z\in M$ the induced morphism
$$Map_{M}(z,x) \longrightarrow Map_{M}(z,y)$$
induces an injection on $\pi_{0}$ and isomorphisms on all
$\pi_{i}$ for $i>0$ (for all base points). This is also equivalent to say that
the natural morphism
$$x \longrightarrow x\times^{h}_{y}x$$
is an isomorphism in $Ho(M)$.
The following lemma will be used
implicitly in the sequel.

\begin{lem}\label{lmono}
A morphism $f : C \longrightarrow D$ in $dg-Cat$ is a homotopy monomorphism if and
only if it is quasi-fully faithful.
\end{lem}

\textit{Proof:} We can of course suppose that the morphism
$f$ is a fibration in $dg-Cat$. Then, $f$ is a homotopy monomorphism
if and only if the induced morphism
$$\Delta : C \longrightarrow C\times_{D}C$$
is a quasi-equivalence.

Let us first assume that $f$ is quasi-fully faithful.
For any $x$ and $y$ in $C$ the induced morphism by $\Delta$ is
the diagonal of $C(x,y)$
$$\Delta(x,y) : C(x,y) \longrightarrow C(x,y)\times_{D(f(x),f(y))}C(x,y).$$
As $f$ is a fibration, the morphism $C(x,y) \longrightarrow D(f(x),f(y))$
is a trivial fibration, and thus the morphism $\Delta(x,y)$
is a quasi-isomorphism.
This shows that $\Delta$ is quasi-fully faithful.
Now, let $t$ be an object in $C\times_{D}C$, corresponding to
two points $x$ and $y$ in $C$ such that $f(x)=f(y)$. We consider the
identity morphism $f(x) \rightarrow f(y)$ in $[D]$. As $[C] \rightarrow [D]$ is
fully faithful, the identity can be lifted to an isomorphism in $[C]$
$u : x \rightarrow y$. Furthermore, as $C(x,y) \longrightarrow D(f(x),f(y))$
is a fibration, the morphism $u$ can be represented by
a zero cycle $u\in Z^{0}(C(x,y))$ whose image by $f$ is the identity.
This implies that the point $t$ is isomorphic in $[C\times_{D}C]$ to
the image of the point $x\in C$ by $\Delta$, and thus
that $\Delta$ is quasi-essentially surjective. We have shown
that $\Delta$ is a quasi-equivalence and therefore that $f$
is a homotopy monomorphism.

Conversely, let us assume that $f$ is a homotopy monomorphism. Then,
for any $x$ and $y$ in $C$ the natural morphism
$$C(x,y) \longrightarrow C(x,y)\times_{D(f(x),f(y))}C(x,y)$$
is a quasi-isomorphism, and thus the morphism $C(x,y) \longrightarrow D(f(x),f(y))$
is a homotopy monomorphism in $C(k)$. As $C(k)$ is a stable model category
(see \cite[\S 7]{ho})this clearly implies that $C(x,y) \longrightarrow D(f(x),f(y))$
is in fact a quasi-isomorphism. \hfill $\Box$ \\

\begin{cor}\label{clmono}
Let $C \longrightarrow D$ be a quasi-fully faithful morphism in $dg-Cat$
and $B$ be any dg-category. Then, the induced morphism
$$Map(B,C) \longrightarrow Map(B,D)$$
induces an injection on $\pi_{0}$ and an isomorphism
on $\pi_{i}$ for $i>0$. Furthermore, the image of
$$\pi_{0}(Map(B,C))=[B,C] \longrightarrow
[B,D]=\pi_{0}(Map(B,D))$$ consists of all morphism such that the
induced functor $[B] \rightarrow [D]$ factors through
the essential image of $[C] \rightarrow [D]$.
\end{cor}

\textit{Proof:} Only the last statement requires a proof.
For this we can of course assume that $B$ is cofibrant. Furthermore, one can
replace $C$ by its quasi-essential image in $D$. The statement is then clear
by the description of $[B,C]$ and $[B,D]$ as
homotopy classes of morphisms between $B$ and $C$ or $D$. \hfill $\Box$ \\

\section{Modules over dg-categories}

Let $C \in dg-Cat_{\mathbb{U}}$
be a fixed $\mathbb{U}$-small $dg$-category. Recall that
a $\mathbb{U}$-small $C$-dg-module $F$ (or simply
a $C$-module) consists of the following data.

\begin{itemize}
\item For any object $x\in C$ a complex $F(x)\in C(k)_{\mathbb{U}}$.
\item For any two objects
$x$ and $y$ in $C$, a morphism of complexes
$$ C(x,y)\otimes F(x) \longrightarrow F(y),$$
satisfying the usual associativity and unit conditions.
\end{itemize}

Note that a $C$-module is nothing else than
a morphism of dg-categories
$F : C \longrightarrow C(k)$, where
$C(k)$ is a dg-category in the obvious way, or equivalently
as a $C(k)$-enriched functor from $C$ to $C(k)$.
For two $C$-dg-modules $F$ and $G$, a morphism
from $F$ to $G$ is simply the data of
morphisms $f_{x} : F(x) \longrightarrow G(x)$
commuting with the structure morphisms. This
is nothing else than a $C(k)$-enriched natural transformation
between the corresponding
$C(k)$-enriched functors. The $\mathbb{U}$-small $C$-modules
and morphisms between them form
a category, denoted by $C-Mod_{\mathbb{U}}$.
Once again, when the universe $\mathbb{U}$ is irrelevant
we will simply write $C-Mod$ for $C-Mod_{\mathbb{U}}$. \\

Let $z\in C$ be an object in $C$. One defines
a $C$-module
$\underline{h}^{z}\in C-Mod$, by the formula
$\underline{h}^{z}(x):=C(z,x)$, and with
structure morphisms
$$C(z,x)\otimes C(x,y) \longrightarrow C(z,y)$$
being the composition in $C$.

\begin{df}\label{d2}
Let $C\in dg-Cat$ and $f : F \longrightarrow G$ be a morphism
of $C$-modules. The morphism $f$ is
an \emph{equivalence} (resp. a \emph{fibration})
if for any $x\in C$ the morphism
$$f_{x} : F(x) \longrightarrow G(x)$$
is an equivalence (resp. a fibration) in $C(k)$.
\end{df}

We recall that as $C(k)$ is cofibrantly generated,
the above definition endows $C-Mod$ with a structure of a
cofibrantly generated model category (see for example
\cite[\S 11]{hi}). The natural $C(k)$-enrichment of
$C-Mod$ endows furthermore $C-Mod$ with a structure
of a $C(k)$-model category in the sense of \cite[4.2.18]{ho}.
The $C(k)$-enriched $Hom$'s of the category
$C-Mod$ will be denoted by $\underline{Hom}$, and its
derived version by
$$\mathbb{R}\underline{Hom} : Ho(C-Mod)^{op}\times Ho(C-Mod) \longrightarrow
Ho(C(k)).$$

The notion of modules over dg-categories
has the following natural  generalization.
Let $M$ be a $C(k)_{\mathbb{U}}$-model category in the sense
of \cite[4.2.18]{ho}, and let us suppose that it
is $\mathbb{U}$-cofibrantly generated in the
sense of \cite[Appendix A]{hagI}. Then, for a
$\mathbb{U}$-small dg-category
$C$ one has a category of $C(k)$-enriched
functors $M^{C}$ from $C$ to $M$.
Furthermore, it can be endowed with a structure
of a $\mathbb{U}$-cofibrantly generated model
category for which equivalences and fibrations are
defined levelwise in $M$ (see e.g. \cite[11.6]{hi}).
The category
$M^{C}$ has itself a natural $C(k)$-enrichment induced
from the one on $M$,
making it into a $C(k)$-model category.
When $M=C(k)_{\mathbb{U}}$ itself, the model category
$M^{C}$ can be identified with $C-Mod_{\mathbb{U}}$.

Let $f : C \longrightarrow D$ be a morphism in $dg-Cat$.
Composing with $f$ gives a restriction functor
$$f^{*} : M^{D} \longrightarrow M^{C}.$$
This functor has a left adjoint
$$f_{!} : M^{C} \longrightarrow M^{D}.$$
The adjunction $(f_{!},f^{*})$ is clearly a Quillen adjunction, compatible
with the $C(k)$-enrichment.

\begin{prop}\label{p1}
Let $f : C \longrightarrow D$ be a quasi-equivalence
between $\mathbb{U}$-small dg-categories. Let
$M$ be a $\mathbb{U}$-cofibrantly generated
$C(k)$-model category, such that the
domain and codomain of a set of generating
cofibrations are cofibrant objects in $M$. We assume that
one of the following conditions is satisfied.
\begin{enumerate}
\item For
any cofibrant object $A\in M$, and any
quasi-isomorphism $X \longrightarrow Y$ in $C(k)$, the
induced morphism
$$X\otimes A \longrightarrow Y\otimes A$$
is an equivalence in $M$.
\item All the complexes of morphisms of $C$ and
$D$ are cofibrant objects in $C(k)$.
\end{enumerate}
Then
the Quillen adjunction $(f_{!},f^{*})$ is
a Quillen equivalence.
\end{prop}

\textit{Proof:} The functor $f^{*}$ clearly preserves
equivalences. Furthermore, as $f$ is quasi-essentially
surjective,
the functor $f^{*} : Ho(M^{D}) \longrightarrow Ho(M^{C})$
is easily seen to be conservative. Therefore, one is reduced to check that
the adjunction morphism $Id\Rightarrow  f^{*}\mathbb{L}f_{!} $
is an isomorphism.

For $x\in C$, and $A\in M$, one writes
$\underline{h}^{x}\otimes A \in M^{C}$ for the object
defined by
$$\begin{array}{cccc}
\underline{h}^{x}\otimes A & C & \longrightarrow & M \\
 & y & \mapsto & C(x,y)\otimes A.
\end{array}$$
The model category $M^{C}$ is itself
cofibrantly generated, and a set of
generating cofibration can be chosen
to consist of morphisms of the form
$$\underline{h}^{x}\otimes A \longrightarrow \underline{h}^{x}\otimes B$$
for some generating cofibration $A\longrightarrow B$ in $M$.
By assumption on $M$, any object $F\in Ho(M^{C})$ can thus
be written as a
homotopy colimit of objects of the form
$\underline{h}^{x}\otimes A$, for
certain cofibrant $A\in M$, and certain $x\in C$. As
the two functors $f^{*}$ and $\mathbb{L}f_{!}$ commute
with homotopy colimits it is then enough to show that
the natural morphism
$$\underline{h}^{x}\otimes A
\longrightarrow f^{*}\mathbb{L}f_{!}
(\underline{h}^{x}\otimes A)
$$
is an isomorphism in $Ho(M^{C})$. By adjunction, one clearly has
$\mathbb{L}f_{!}(\underline{h}^{x}\otimes A)\simeq
\underline{h}^{f(x)}\otimes A$.
Therefore, the adjunction morphism
$$\underline{h}^{x}\otimes A
\longrightarrow f^{*}\mathbb{L}f_{!}
(\underline{h}^{x}\otimes A)\simeq f^{*}(\underline{h}^{f(x)}\otimes A)
$$
evaluated at $y\in C$ is the morphism
$$f_{x,y}\otimes Id_{A} : C(x,y) \otimes A \longrightarrow
D(f(x),f(y))\otimes A.$$
The fact that this is an isomorphism in $Ho(M)$
follows from the fact that
$f$ is quasi-fully faithful, one of our hypothesis $(1)$ and $(2)$, and the fact that
$M$ is a $C(k)$-model category.  \hfill $\Box$ \\

Another important property of the model category
$M^{C}$ is the following.

\begin{prop}\label{p2}
Let $C$ be a $\mathbb{U}$-small dg-category
with cofibrant complexes of morphisms
(i.e. $C(x,y)$ is cofibrant in $C(k)$ for all
$x$ and $y$),
and $M$ be a  $\mathbb{U}$-cofibrantly generated
$C(k)$-model category. Then,
for any $x\in C$ the evaluation functor
$$\begin{array}{cccc}
x^{*} : & M^{C} & \longrightarrow & M \\
 & F & \mapsto & F(x)
\end{array}$$
preserves fibrations, cofibrations and equivalences.
\end{prop}

\textit{Proof:} For fibrations and equivalences this
is clear by definition. The
functor $x^{*}$ commutes with colimits, and thus
by a small object argument one is
reduced to show that $x^{*}$ sends
generating cofibrations to cofibrations.
One
knows that the generating set of cofibrations
in $M^{C}$ can be chosen to consist of morphisms
of the form $\underline{h}^{z}\otimes A \longrightarrow
\underline{h}^{z}\otimes B$ for some cofibration
$A \longrightarrow B$ in $M$. The image by $x^{*}$
of such a morphism is
$$C(z,x)\otimes A  \longrightarrow
 C(z,x)\otimes B.$$
As by assumption $C(z,x)$ is a cofibrant object in $C(k)$, one
sees that this morphism is a cofibration in $M$. \hfill $\Box$ \\

Two important cases of application of proposition \ref{p2} is when
$C$ itself is a cofibrant dg-category (see Prop. \ref{p0}), or
when $k$ is a field. \\

\begin{cor}\label{cp2}
The conclusion of Prop. \ref{p1} is satisfied
when $M$ is of the form $D-Mod_{\mathbb{U}}$, for
a $\mathbb{U}$-small dg-category $D$ with
cofibrant complexes of morphisms (in particular for
$M=C(k)$).
\end{cor}

\textit{Proof:} This follows easily from Prop. \ref{p2} and the
fact that $C(k)$ itself satisfies the hypothesis $(1)$
of Prop. \ref{p1}. \hfill $\Box$ \\

Let $\mathbb{U} \in \mathbb{V}$ be two universes. Let $M$ be
a $C(k)_{\mathbb{U}}$-model category which is supposed to be
furthermore $\mathbb{V}$-small.
We define a $\mathbb{V}$-small dg-category
$Int(M)$ in the following way\footnote{The notation \emph{Int} is taken from
\cite{hs}. As far as I understand it stands for \emph{internal}.}. The set of objects
of $Int(M)$ is the set of fibrant and cofibrant objects
in $M$. For two such objects $F$ and $E$
one sets
$$Int(M)(E,F):=\underline{Hom}(E,F)\in C(k)_{\mathbb{U}},$$
where $\underline{Hom}(E,F)$ is the $C(k)$-valued
Hom of the category $M$.
The dg-category $Int(M)$ is of course
only $\mathbb{V}$-small as its sets of objects is only $\mathbb{V}$-small.
However, for any $E$ and $F$ in $Int(M)$
the complex $Int(M)(E,F)$ is in fact
$\mathbb{U}$-small. \\

The following is a general fact about
$C(k)$-enriched model categories.

\begin{prop}\label{p3}
There exists a natural equivalence of categories
$$[Int(M)]\simeq Ho(M).$$
\end{prop}

\textit{Proof:} This follows from the formula
$$H^{0}(\mathbb{R}\underline{Hom}(X,Y))\simeq [k,\mathbb{R}\underline{Hom}(X,Y)]_{C(k)}\simeq
[X,Y]_{M},$$
for two objects $X$ and $Y$ in $M$. \hfill $\Box$ \\

For $x\in C$, the object $\underline{h}^{x} \in C-Mod_{\mathbb{U}}$ is
cofibrant and fibrant, and therefore
the construction $x \mapsto
\underline{h}^{x}$, provides a morphism of dg-categories
$$\underline{h}^{-} : C^{op} \longrightarrow Int(C-Mod_{\mathbb{U}}),$$
where $C^{op}$ is the opposite dg-category of $C$ ($C^{op}$ has the same set of objects than $C$ and
$C^{op}(x,y):=C(y,x)$).
The morphism $\underline{h}^{-}$ can also be written dually as
$$\underline{h}_{-} : C \longrightarrow
Int(C^{op}-Mod_{\mathbb{U}}).$$
The dg-functor $\underline{h}^{-}$ will be considered as
a morphism in $dg-Cat_{\mathbb{V}}$, and is clearly
quasi-fully faithful by an application of the $C(k)$-enriched
Yoneda lemma.

\begin{df}\label{d3}
\begin{enumerate}
\item
Let $C\in dg-Cat_{\mathbb{U}}$, and
$F\in C^{op}-Mod_{\mathbb{U}}$ be a $C^{op}$-module.
The object $F$ is called \emph{representable} (resp.
\emph{quasi-representable}) if it is isomorphic
in $C^{op}-Mod_{\mathbb{U}}$ (resp.
in $Ho(C^{op}-Mod_{\mathbb{U}})$) to $\underline{h}_{x}$
for some object $x\in C$.
\item Dually, let $C\in dg-Cat_{\mathbb{U}}$, and
$F\in C-Mod_{\mathbb{U}}$ be a $C$-module.
The object $F$ is called \emph{corepresentable} (resp.
\emph{quasi-corepresentable}) if it is isomorphic
in $C-Mod_{\mathbb{U}}$ (resp.
in $Ho(C-Mod_{\mathbb{U}})$) to $\underline{h}^{x}$
for some object $x\in C$.
\end{enumerate}
\end{df}

As the morphism $\underline{h}_{-}$ is quasi-fully faithful,
it induces a quasi-equivalence between $C$ and the
full dg-category of $Int(C^{op}-Mod_{\mathbb{U}})$
consisting of quasi-representable objects. This quasi-equivalence is
a morphism in $dg-Cat_{\mathbb{V}}$. \\

\section{Mapping spaces and bi-modules}

Let $C$ and $D$ be two objects in $dg-Cat$. One has
a tensor product $C\otimes D \in dg-Cat$ defined in the following way.
The set of objects of $C\otimes D$ is $Ob(C)\times Ob(D)$, and
for $(x,y)$ and $(x',y')$ two objects in $Ob(C\otimes D)$ one sets
$$(C\otimes D)((x,y),(x',y')):=C(x,y)\otimes D(x',y').$$
Composition in $C\otimes D$ is given by the obvious
formula.
This defines a symmetric monoidal structure on
$dg-Cat$, which is easily seen to be
closed. The unit of this structure will be
denoted by $\mathbf{1}$, and is the dg-category with
a unique object and $k$ as its endomorphism ring.

The model category $dg-Cat$ together with
the symmetric monoidal structure $-\otimes -$ is \emph{not}
a symmetric monoidal model category, as
the tensor product of two cofibrant objects in
$dg-Cat$ is not cofibrant in general. A direct consequence
of this fact is that the internal Hom object
between cofibrant-fibrant objects in $dg-Cat$
can not be invariant by quasi-equivalences, and thus does not
provide internal Hom's for the homotopy categories
$Ho(dg-Cat)$. This fact is the main difficulty in
computing the mapping spaces in $dg-Cat$, as the naive
approach simply does not work.

However, it is true that the monoidal structure
$\otimes$ on $dg-Cat$ is closed, and that
$dg-Cat$ has corresponding internal Hom objects
$C^{D}$ satisfying the usual adjunction rule
$$Hom_{dg-Cat}(A\otimes B,C)\simeq Hom(A,C^{B}).$$
This gives a natural equivalence of categories
$$M^{C\otimes D}\simeq (M^{C})^{D}$$
for any $C(k)$-enriched category $M$. Furthermore, when
$M$ is a $\mathbb{U}$-cofibrantly generated model category,
this last equivalence is compatible with the
model structures on both sides. \\

The functor $-\otimes -$ can be derived into a functor
$$-\otimes^{\mathbb{L}} - : dg-Cat\times dg-Cat \longrightarrow dg-Cat$$
defined by the formula
$$C\otimes^{\mathbb{L}}D:=
Q(C)\otimes D$$
where $Q$ is a cofibrant replacement in $dg-Cat$ which acts by
the identity on the sets of objects. Clearly,
the functor $-\otimes^{\mathbb{L}} -$ preserves quasi-equivalences and
passes through the homotopy categories
$$-\otimes^{\mathbb{L}} - : Ho(dg-Cat)\times Ho(dg-Cat) \longrightarrow
Ho(dg-Cat).$$
Note that when $C$ is cofibrant, one has
a natural quasi-equivalence $C\otimes^{\mathbb{L}}D \longrightarrow C\otimes D$.

We now consider $(C\otimes D^{op})-Mod$, the category of
$(C\otimes D^{op})$-modules. For any
object $x\in C$, there exists a natural
morphism of dg-categories $D^{op} \longrightarrow  (C\otimes D^{op})$ sending
$y\in D$ to the object $(x,y)$, and
$$D^{op}(y,z) \longrightarrow (C\otimes D^{op})((x,y),(x,z))=
C(x,x)\otimes D^{op}(y,z)$$
being the tensor product of the unit $k \longrightarrow C(x,x)$ and
the identity on $D^{op}(y,z)$. As $C$ and $Q(C)$ has the same
set of objects, one sees that for any $x\in C$ one also gets
a natural morphism of dg-categories
$$i_{x} : D^{op} \longrightarrow Q(C)\otimes D^{op}=C\otimes^{\mathbb{L}}D^{op}.$$

\begin{df}\label{d4}
Let $C$ and $D$ be two dg-categories.
An object $F \in (C\otimes^{\mathbb{L}}D^{op})-Mod$ is called
\emph{right quasi-representable}, if for any $x\in C$, the $D^{op}$-module
$i_{x}^{*}(F)\in D^{op}-Mod$ is quasi-representable in the sense
of Def. \ref{d3}.
\end{df}

We now let $\mathbb{U}\in \mathbb{V}$ be two universes, and
let $C$ and $D$ be two $\mathbb{U}$-small dg-categories.
Let
$\Gamma^{*}$ be a co-simplicial resolution functor in $dg-Cat_{\mathbb{U}}$
in the sense of \cite[\S 16.1]{hi}. Recall that
$\Gamma^{*}$ is a functor from $dg-Cat_{\mathbb{U}}$ to
$dg-Cat^{\Delta}_{\mathbb{U}}$, equipped with a natural
augmentation $\Gamma^{0} \longrightarrow Id$, and such the
following two conditions are satisfied.
\begin{itemize}
\item For any $n$, and any $C\in dg-Cat_{\mathbb{U}}$ the morphism
$\Gamma^{n}(C) \rightarrow C$
is a quasi-equivalence.
\item For any $C\in dg-Cat_{\mathbb{U}}$, the object
$\Gamma^{*}(C) \in dg-Cat_{\mathbb{U}}^{\Delta}$ is cofibrant
for the Reedy model structure.
\item The morphism $\Gamma^{0}(C) \longrightarrow C$
is equal to $Q(C) \longrightarrow C$.
\end{itemize}

The left mapping space between $C$ and $D$ is by definition the
$\mathbb{U}$-small simplicial set
$$\begin{array}{cccc}
Map^{l}(C,D):=Hom(\Gamma^{*}(C),D) : & \Delta^{op} & \longrightarrow & Set_{\mathbb{U}} \\
 & [n] & \mapsto & Hom(\Gamma^{n}(C),D).
\end{array}$$
Note that the mapping space $Map^{l}(C,D)$ defined above has
the correct homotopy type as all objects are fibrant
in $dg-Cat_{\mathbb{U}}$.

For any $[n] \in \Delta$, one considers the
(non-full) sub-category $\mathcal{M}(\Gamma^{n}(C),D)$
of $(\Gamma^{n}(C)\otimes D^{op})-Mod_{\mathbb{U}}$
defined in the following way. The objects
of $\mathcal{M}(\Gamma^{n}(C),D)$ are the
$(\Gamma^{n}(C)\otimes D^{op})$-modules $F$ such that
$F$ is right quasi-representable, and for any
$x\in \Gamma^{n}(C)$ the $D^{op}$-module $F(x,-)$
is cofibrant in $D^{op}-Mod_{\mathbb{U}}$. The morphisms
in $\mathcal{M}(\Gamma^{n}(C),D)$ are simply the
equivalences in $(\Gamma^{n}(C)\otimes D^{op})-Mod_{\mathbb{U}}$.
The nerve of the category $\mathcal{M}(\Gamma^{n}(C),D)$ gives
a $\mathbb{V}$-small simplicial set
$N(\mathcal{M}(\Gamma^{n}(C),D))$. For $[n] \rightarrow [m]$
a morphism in $\Delta$, one has a natural morphism
of dg-categories $\Gamma^{n}(C)\otimes D^{op} \longrightarrow
\Gamma^{m}(C)\otimes D^{op}$, and thus a well defined
morphism of simplicial sets
$$N(\mathcal{M}(\Gamma^{m}(C),D)) \longrightarrow N(\mathcal{M}(\Gamma^{n}(C),D))$$
obtained by pulling back the modules from
$\Gamma^{m}(C)\otimes D^{op}$ to $\Gamma^{n}(C)\otimes D^{op}$. This
defines a functor
$$\begin{array}{cccc}
N(\mathcal{M}(\Gamma^{*}(C),D)) : & \Delta^{op} & \longrightarrow &
SSet_{\mathbb{V}} \\
 & [n] & \mapsto & N(\mathcal{M}(\Gamma^{n}(C),D)).
\end{array}$$

The set of zero simplices in $N(\mathcal{M}(\Gamma^{n}(C),D))$
is the set of all objects in the category
$\mathcal{M}(\Gamma^{n}(C),D)$. Therefore,
one defines a natural morphism of sets
$$Hom(\Gamma^{n}(C),D) \longrightarrow
N(\mathcal{M}(\Gamma^{n}(C),D))_{0}$$
by sending a morphism of dg-categories $f : \Gamma^{n}(C) \longrightarrow D$,
to the $(\Gamma^{n}(C)\otimes D^{op})$-module $\phi(f)$
defined by $\phi(f)(x,y):=D(y,f(x))$ and the natural
transition morphisms. Note that $\phi(f)$ belongs to
the sub-category $\mathcal{M}(\Gamma^{n}(C),D)$ as
for any $x\in \Gamma^{n}(C)$ the $D^{op}$-module
$\phi(f)(x,-)=\underline{h}_{f(x)}$
is representable and thus
quasi-representable and cofibrant.
By adjunction, this morphism of sets can also
be considered as a morphism of simplicial sets
$$\phi : Hom(\Gamma^{n}(C),D) \longrightarrow
N(\mathcal{M}(\Gamma^{n}(C),D)),$$
where the set $Hom(\Gamma^{n}(C),D)$ is considered as a constant
simplicial set.
This construction is clearly functorial
in $n$, and gives a well defined morphism of
bi-simplicial sets
$$\phi : Hom(\Gamma^{*}(C),D) \longrightarrow N(\mathcal{M}(\Gamma^{*}(C),D)).$$
Passing to the diagonal one gets a morphism in
$SSet_{\mathbb{V}}$
$$\phi : Map^{l}(C,D) \longrightarrow
d(N(\mathcal{M}(\Gamma^{*}(C),D))).$$
Finally, the diagonal $d(N(\mathcal{M}(\Gamma^{*}(C),D)))$ receives a natural
morphism
$$\psi : N(\mathcal{M}(\Gamma^{0}(C),D))=N(\mathcal{M}(Q(C),D)) \longrightarrow
d(N(\mathcal{M}(\Gamma^{*}(C),D))).$$
Clearly, the diagram of simplicial sets
$$\xymatrix{
Map^{l}(C,D) \ar[r]^-{\phi} & d(N(\mathcal{M}(\Gamma^{*}(C),D))) &
\ar[l]_-{\psi} N(\mathcal{M}(Q(C),D))}$$
is functorial in $C$. \\

The main theorem of this work is the following.

\begin{thm}\label{t1}
The two morphisms in $SSet_{\mathbb{V}}$
$$\xymatrix{
Map^{l}(C,D) \ar[r]^-{\phi} & d(N(\mathcal{M}(\Gamma^{*}(C),D))) &
\ar[l]_-{\psi} N(\mathcal{M}(Q(C),D))}$$
are weak equivalences.
\end{thm}

\textit{Proof:} For any $n$, the morphism
$\Gamma^{n}(C)\otimes D^{op} \longrightarrow Q(C)\otimes D^{op}$
is a quasi-equivalence of dg-categories. Therefore,
Prop. \ref{cp2} implies that the pull-back functor
$$(Q(C)\otimes D^{op})-Mod \longrightarrow
(\Gamma^{n}(C)\otimes D^{op})-Mod$$
is the right adjoint of a Quillen equivalence. As these functors
obviously preserve the notion of being right quasi-representable,
one finds that the induced morphism
$$N(\mathcal{M}(Q(C),D)) \longrightarrow
N(\mathcal{M}(\Gamma^{n}(C),D))$$
is a weak equivalence. This clearly implies that
the morphism $\psi$ is a weak equivalence.

It remains to show that the morphism $\phi$ is also
a weak equivalence. For this, we start by proving that
it induces an isomorphism on connected components.

\begin{lem}\label{lt1}
The induced morphism
$$\pi_{0}(\phi) : [C,D]\simeq \pi_{0}(Map^{l}(C,D))\longrightarrow
\pi_{0}(d(N(\mathcal{M}(\Gamma^{*}(C),D))))$$
is an isomorphism.
\end{lem}

\textit{Proof:} First of all, replacing $C$ by
$Q(C)$ one can suppose that $Q(C)=C$
(one can do this because of Prop. \ref{cp2}).
One then has $\pi_{0}(Map^{l}(C,D))\simeq [C,D]$, and
$\pi_{0}(d(N(\mathcal{M}(\Gamma^{*}(C),D))))\simeq
\pi_{0}(N(\mathcal{M}(C,D)))$ is the set of
isomorphism classes in $Ho((C\otimes D^{op})-Mod)^{rqr}$,
the full sub-category of $Ho((C\otimes D^{op})-Mod)$
consisting of all right quasi-representable objects. The morphism
$\phi$ naturally gives a morphism
$$\phi : [C,D] \longrightarrow Iso(Ho((C\otimes D^{op})-Mod)^{rqr})$$
which can be described as follows. For
any $f\in [C,D]$, represented by $f : C \longrightarrow D$
in $Ho(dg-Cat)$, $\phi(f)$ is the $C\otimes D^{op}$-module
defined by $\phi(f)(x,y):=D(y,f(x))$.

\begin{sublem}\label{sl2}
With the same notations as above,
let
$M$ be a $\mathbb{U}$-cofibrantly generated
$C(k)_{\mathbb{U}}$-model category, which is
furthermore $\mathbb{V}$-small. Let
$Iso(Ho(M^{C}))$ be the set of isomorphism classes
of objects in $Ho(M^{C})$.
Then, the natural morphism
$$Hom(C,Int(M)) \longrightarrow Iso(Ho(M^{C}))$$
is surjective.
\end{sublem}

\textit{Proof of sub-lemma \ref{sl2}:} Of course, the morphism
$$Hom(C,Int(M)) \longrightarrow Iso(Ho(M^{C}))$$
sends a morphism of dg-categories $C \longrightarrow Int(M)$
to the corresponding object in $M^{C}$.
Let $F\in Ho(M^{C})$ be
a any cofibrant and fibrant object.
This object is given by
a $C(k)$-enriched functor $F : C \longrightarrow M$.
Furthermore, as $F$ is fibrant and cofibrant,
Prop. \ref{p2} tells us that $F(x)$ is fibrant and cofibrant in $M$
for any $x\in C$.
The object $F$ can therefore be naturally considered as a morphism
of $\mathbb{V}$-small dg-categories
$$F : C \longrightarrow Int(M),$$
which gives an element in $Hom(C,Int(M))$ sent to $F$ by the map of the lemma.
\hfill $\Box$ \\

Let us now prove that the morphism $\phi$ is surjective
on connected component.
For this, let $F\in Ho((C\otimes D^{op})-Mod_{\mathbb{U}})$ be
a right quasi-representable object. One needs to show that
$F$ is isomorphic to some $\phi(f)$ for some
morphism of dg-categories $f : C \longrightarrow D$.
Sub-lemma \ref{sl2} applied to $M=D^{op}-Mod_{\mathbb{U}}$ implies
that $F$ corresponds to a morphism
of $\mathbb{V}$-small dg-categories
$$F : C \longrightarrow Int(D^{op}-Mod_{\mathbb{U}})^{qr},$$
where $Int(D^{op}-Mod)^{qr}$ is the full sub-dg-category
of $Int(D^{op}-Mod_{\mathbb{U}})$ consisting of all
quasi-representable objects.

One has a diagram
in $dg-Cat_{\mathbb{V}}$
$$\xymatrix{
C \ar[r]^-{F} & Int(D^{op}-Mod_{\mathbb{U}})^{qr} \\
 & D \ar[u]_-{\underline{h}}.}$$
As the morphism $\underline{h}$ is a quasi-equivalence, and
as $C$ is cofibrant, one finds a morphism of dg-categories
$f : C \longrightarrow D$, such that the two morphisms
$$F : C \longrightarrow Int(D^{op}-Mod_{\mathbb{U}})^{qr} \qquad
\underline{h}_{f(-)}=\phi(f) : C \longrightarrow Int(D^{op}-Mod_{\mathbb{U}})^{qr}$$
are homotopic in $dg-Cat_{\mathbb{V}}$. Let
$$\xymatrix{
C \ar[d]_-{i_{0}}\ar[rd]^-{F} & \\
C' \ar[r]^-{H} & Int(D^{op}-Mod_{\mathbb{U}})^{qr} \\
C\ar[u]^-{i_{1}} \ar[ru]_-{\phi(f)} & }$$
be a homotopy in $dg-Cat_{\mathbb{V}}$. Note that
$C'$ is a cylinder object for $C$, and thus can be
chosen to be cofibrant and $\mathbb{U}$-small. We let
$p : C' \longrightarrow C$ the natural projection, such that
$p\circ i_{0}=p\circ i_{1}=Id$. This diagram gives rise to an
equivalence of categories (by Prop. \ref{cp2})
$$i_{0}^{*}\simeq i_{1}^{*}\simeq (p^{*})^{-1} : Ho((C'\otimes D^{op})-Mod_{\mathbb{U}}) \longrightarrow
Ho((C\otimes D^{op})-Mod_{\mathbb{U}}).$$
Furthermore, one has
$$F\simeq i_{0}^{*}(H)\simeq i_{1}^{*}(H)\simeq \phi(f).$$
This shows that the two $C\otimes D^{op}$-modules
$F$ and $\phi(f)$ are
isomorphic in $Ho((C\otimes D^{op})-Mod_{\mathbb{U}})$, or in other
words that $\phi(f)=F$ in $Iso(Ho((C\otimes D^{op})-Mod_{\mathbb{U}}))$.
This finishes the proof of the surjectivity part of the
lemma \ref{lt1}. \\

Let us now prove that $\phi$ is injective. For this, let
$f,g : C \longrightarrow D$ be two morphisms of dg-categories, such that
the two $(C\otimes D^{op})$-modules $\phi(f)$ and $\phi(g)$ are
isomorphic in $Ho((C\otimes D^{op})-Mod_{\mathbb{U}}))$. Composing
$f$ and $g$ with
$$\underline{h} : D \longrightarrow Int(D^{op}-Mod_{\mathbb{U}})$$
one gets two new morphisms of dg-categories
$$f',g' : C \longrightarrow Int(D^{op}-Mod_{\mathbb{U}}).$$
Using that $\underline{h}$ is quasi-fully faithful
Cor. \ref{clmono}
implies that
if $f'$ and $g'$ are homotopic morphisms in $dg-Cat_{\mathbb{V}}$, then
$f$ and $g$ are equal as morphisms in $Ho(dg-Cat_{\mathbb{V}})$.
As the inclusion $Ho(dg-Cat_{\mathbb{U}}) \longrightarrow
Ho(dg-Cat_{\mathbb{V}})$ is fully faithful
(see remark after Def. \ref{d1}), we see it is enough to show that
$f'$ and $g'$ are homotopic in $dg-Cat_{\mathbb{V}}$.

\begin{sublem}\label{sl1}
Let $M$ be a $C(k)_{\mathbb{U}}$-model category which is
$\mathbb{U}$-cofibrantly generated and $\mathbb{V}$-small.
Let $u$ and $v$ be two morphisms in $dg-Cat_{\mathbb{V}}$
$$u,v : C \longrightarrow Int(M)$$
such that the corresponding objects
in $Ho(M^{C})$ are isomorphic. Then
$u$ and $v$ are homotopic as morphisms in $dg-Cat_{\mathbb{V}}$.
\end{sublem}

\textit{Proof of sub-lemma \ref{sl1}:} First of all, any isomorphism in
$Ho(M^{C})$ between levelwise cofibrant and fibrant objects can be represented as a string of
trivial cofibrations and trivial fibrations between levelwise cofibrant and fibrant objects.
Therefore, sub-lemma \ref{sl2} shows that one is
reduced
to the case where there exists
an equivalence $\alpha : u \longrightarrow v$ in
$M^{C}$ which is either a fibration or a cofibration.

Let us start with the case where $\alpha$ is a cofibration
in $M^{C}$. The morphism $\alpha$ can also be considered
as an object in $(M^{I})^{C}$, where
$I$ is the category with two
objects $0$ and $1$ and a unique morphism
$0\rightarrow 1$. The category
$M^{I}$, which is the category
of morphisms in $M$, is endowed with its projective model structure, for
which fibrations and equivalences are defined on the underlying
objects in $M$. As the morphism $\alpha$ is a cofibration in
$M^{C}$, we see that for $x\in C$ the
corresponding morphism $\alpha_{x} : u(x) \rightarrow v(x)$
is a cofibration in $M$, and thus is a cofibrant (and fibrant)
object in $M^{I}$ because of proposition \ref{p2}. This implies that $\alpha$ gives rise to
a morphism of dg-categories
$$\alpha : C \longrightarrow Int(M^{I}).$$
Now, let $Int(M) \longrightarrow Int(M^{I})$ be the natural
inclusion morphism, sending a cofibrant and fibrant object in $M$ to the
identity morphism. This a morphism in $dg-Cat_{\mathbb{V}}$ which
is easily seen to be quasi-fully faithful. We let
$C'\subset Int(M^{I})$ be the quasi-essential image of
$Int(M)$ in $Int(M^{I})$. It is easy to check that
$C'$ is the full sub-dg-category of $Int(M^{I})$ consisting of
all objects in $M^{I}$ corresponding to equivalences in $M$.
The morphism $\alpha : C \longrightarrow Int(M^{I})$
thus factors through the sub-dg-category $C'\subset Int(M^{I})$.
The two objects $0$ and $1$ of $I$ give two projections
$$C' \subset Int(M^{I}) \rightrightarrows Int(M),$$
both of them having the natural inclusion $Int(M)
\longrightarrow Int(M^{I})$ as a section.
We have thus constructed a commutative diagram in $dg-Cat_{\mathbb{V}}$
$$\xymatrix{
 & Int(M) \\
C \ar[r]^-{\alpha} \ar[ru]^-{u} \ar[rd]_-{v} & C' \ar[u] \ar[d] \\
 & Int(M)}$$
which provides a homotopy between $u$ and $v$ in $dg-Cat_{\mathbb{V}}$.

For the case where $\alpha$ is a fibration in $M^{C}$, one
uses the same argument, but endowing $M^{I}$ with its
injective model structure, for which equivalences and cofibrations
are defined levelwise. We leave the details to the reader.
\hfill $\Box$ \\

We have finished the proof of sub-lemma \ref{sl1}, which applied
to $M=D^{op}-Mod_{\mathbb{U}}$ finishes the proof of the
injectivity on connected components, and thus of lemma \ref{lt1}.
\hfill $\Box$ \\

In order to finish the proof of the theorem, one uses
the functoriality of the morphisms $\phi$ and $\psi$
with respect to $D$. First of all, the
simplicial set $Map^{l}(C,D)=Hom(\Gamma^{*}(C),D)$ is
obviously functorial in $D$. One thus has a functor
$$\begin{array}{cccc}
Map^{l}(C,-) : & dg-Cat_{\mathbb{U}} & \longrightarrow &
SSet_{\mathbb{V}} \\
 & D & \mapsto & Map^{l}(C,D).
\end{array}$$
The functoriality of $N(\mathcal{M}(C,D))$ in $D$ is slightly
more complicated. Let $u : D \longrightarrow E$ a morphism
in $dg-Cat_{\mathbb{U}}$. One has a functor
$$(Id\otimes u_{!}) : (C\otimes D^{op})-Mod_{\mathbb{U}} \longrightarrow
(C\otimes E^{op})-Mod_{\mathbb{U}}.$$
This functor can also be described as
$$(u_{!})^{C} : (D^{op}-Mod_{\mathbb{U}})^{C} \longrightarrow
(E^{op}-Mod_{\mathbb{U}})^{C},$$
the natural extension of the functor $u_{!} : D^{op}-Mod_{\mathbb{U}}
\longrightarrow E^{op}-Mod_{\mathbb{U}}$.
Clearly, the functor $(u_{!})^{C}$ sends the sub-category
$\mathcal{M}(C,D)$ to the sub-category
$\mathcal{M}(C,E)$ (here one uses that
$(u_{!})^{C}$ preserves equivalences
because the object $F\in \mathcal{M}(C,D)$
are such that $F(x,-)$ is cofibrant in $D^{op}-Mod_{\mathbb{U}}$). Unfortunately,
this does not define a
presheaf of categories $\mathcal{M}(C,-)$ on
$dg-Cat_{\mathbb{U}}$, as for two
morphisms
$$\xymatrix{
D \ar[r]^-{u} & E \ar[r]^-{v} & F}$$
of dg-categories one only has a natural isomorphism
$(v\circ u)_{!}\simeq (v_{!})\circ u_{!}$ which in general
is not an identity. However, these natural isomorphisms
makes $D \mapsto \mathcal{M}(C,D)$ into a
lax functor from $dg-Cat_{\mathbb{U}}$ to $Cat_{\mathbb{V}}$.
Using the standard rectification procedure, one can
replace up to a natural equivalence the lax functor
$\mathcal{M}(C,-)$ by a true presheaf of categories
$\mathcal{M}'(C,-)$.
Furthermore, the natural morphism
$$Hom(C,D) \longrightarrow \mathcal{M}(C,D)$$
from the set of morphisms $Hom(C,D)$, considered
as a discrete category, to the
category $\mathcal{M}(C,D)$ clearly gives a
morphism of lax functors
$$Hom(C,-) \longrightarrow \mathcal{M}(C,-).$$
By rectification this also induces a natural morphism
of presheaves of categories
$$Hom(C,-) \longrightarrow \mathcal{M}'(C,-).$$
Passing to the nerve one gets a morphism of functors
from $dg-Cat_{\mathbb{U}}$ to $SSet_{\mathbb{V}}$
$$Hom(C,-) \longrightarrow N(\mathcal{M}'(C,-)).$$
This morphism being functorial in $C$ give a
diagram in $(SSet_{\mathbb{V}})^{dg-Cat_{\mathbb{U}}}$
$$\xymatrix{
Map^{l}(C,-)=Hom(\Gamma^{*}(C),-) \ar[r]^-{\phi'} &
d(N(\mathcal{M}'(\Gamma^{*}(C),-))) & \ar[l]_-{\psi'}
N(\mathcal{M}'(Q(C),-)).}$$
These morphisms, evaluated at an object $D\in dg-Cat_{\mathbb{U}}$
gives a diagram of simplicial sets
$$\xymatrix{
Map^{l}(C,D) \ar[r] &
d(N(\mathcal{M}'(\Gamma^{*}(C),D)) & \ar[l]
N(\mathcal{M}'(Q(C),D)),}$$
weakly equivalent to the diagram
$$\xymatrix{
Map^{l}(C,D) \ar[r] &
d(N(\mathcal{M}(\Gamma^{*}(C),D)) & \ar[l]
N(\mathcal{M}(Q(C),D)).}$$
In order to finish the proof of the theorem it is therefore
enough to show that the two
morphism $\phi'$ and $\psi'$ are weak equivalences
of diagrams of simplicial sets. We already know that
$\psi'$ is a weak equivalence, and thus we obtain a
morphism well defined in $Ho((SSet_{\mathbb{V}})^{dg-Cat_{\mathbb{U}}})$
$$k : (\psi')^{-1}\circ \phi :
Map^{l}(C,-) \longrightarrow N(\mathcal{M}'(Q(C),-)).$$
Using our corollary \ref{cp2} it is easy to see that
the functor $N(\mathcal{M}'(Q(C),-))$ sends quasi-equivalences
to  weak equivalences. Furthermore, the standard
properties of mapping spaces imply that so does the
functor $Map^{l}(C,-)$. 

\begin{sublem}\label{sl3}
Let $k : F \longrightarrow G$ be a morphism in
$(SSet_{\mathbb{V}})^{dg-Cat_{\mathbb{U}}}$. Assume the following
conditions are satisfied.
\begin{enumerate}
\item  Both functors $F$ and $G$ send quasi-equivalences
to weak equivalences.
\item
For any diagram in $dg-Cat_{\mathbb{U}}$
$$\xymatrix{
& C\ar[d]^-{p} \\
D \ar[r] &E}$$
with $p$ a fibration, the commutative diagrams
$$\xymatrix{
F(C\times_{E}D) \ar[r] \ar[d] & F(C) \ar[d] & & & G(C\times_{E}D) \ar[r] \ar[d] &  G(C) \ar[d] \\
F(D) \ar[r] & F(E) & & & G(D) \ar[r] & G(E)}$$
are homotopy cartesian.
\item $F(*)\simeq G(*)\simeq *$, where $*$ is the final
object in $dg-Cat$.
\item For any $C\in dg-Cat_{\mathbb{U}}$ the morphism
$k_{C} : \pi_{0}(F(C)) \longrightarrow \pi_{0}(G(C))$
is an isomorphism.
\end{enumerate}
Then, for any $C\in dg-Cat_{\mathbb{U}}$ the natural
morphism
$$k_{C} : F(C) \longrightarrow G(C)$$
is a weak equivalence.
\end{sublem}

\textit{Proof of sub-lemma \ref{sl3}:} Condition
$(1)$ implies that the induced functors
$$Ho(F),Ho(G) : Ho(dg-Cat_{\mathbb{U}}) \longrightarrow
Ho(SSet_{\mathbb{V}})$$
have natural structures of $Ho(SSet_{\mathbb{U}})$-enriched
functors (see for example \cite[Thm. 2.3.5]{hagI}). In particular, for any
$K\in Ho(SSet_{\mathbb{U}})$, and any $C\in Ho(dg-Cat_{\mathbb{U}})$
one has natural morphisms in $Ho(SSet_{\mathbb{U}})$
$$F(C^{\mathbb{R}K})\longrightarrow
Map(K,F(C)) \qquad G(C^{\mathbb{R}K})\longrightarrow
Map(K,G(C)).$$
Our hypothesis $(2)$ and $(3)$ tells us that
when $K$ is a finite simplicial set,
these morphisms are in fact isomorphisms, as
the object $C^{\mathbb{R}K}$ can be functorially constructed using
successive homotopy products and homotopy fiber products. Therefore, conditions
$(4)$ implies that for any finite $K\in Ho(SSet_{\mathbb{U}})$
and any $C\in dg-Cat_{\mathbb{U}}$, the morphism $k_{C}$
induces an isomorphism
$$k_{C^{\mathbb{R}K}} : \pi_{0}(F(C^{\mathbb{R}K}))\simeq
[K,F(C)] \longrightarrow [K,G(C)]\simeq \pi_{0}(G(C^{\mathbb{R}K})).$$
This of course implies that $F(C) \longrightarrow G(C)$
is a weak equivalence. \hfill $\Box$ \\

In order to finish the proof of theorem \ref{t1} it remains to
show that the two functors $Map^{l}(C,-)$ and
$N(\mathcal{M}'(Q(C),-))$  satisfy the conditions
of sub-lemma \ref{sl3}. The case of
$Map^{l}(C,-)$ is clear by the standard properties
of mapping spaces (see \cite[\S 5.4]{ho} or \cite[\S 17]{hi}).
It only remains to show property
$(2)$ of sub-lemma \ref{sl3} for the functor $N(\mathcal{M}'(Q(C),-))$.

\begin{sublem}\label{sl5}
Let $C$ be a cofibrant $\mathbb{U}$-small dg-category,
and let
$$\xymatrix{
D \ar[r]^-{u} \ar[d]_-{v} & D_{1} \ar[d]^-{p} \\
D_{2} \ar[r]_-{q} & D_{3}}$$
be a cartesian diagram in $dg-Cat_{\mathbb{U}}$
with $p$ a fibration. Then, the
square
$$\xymatrix{
N(\mathcal{M}'(C,D)) \ar[r] \ar[d] & N(\mathcal{M}'(C,D_{1})) \ar[d] \\
N(\mathcal{M}'(C,D_{2})) \ar[r] & N(\mathcal{M}'(C,D_{3}))}$$
is homotopy cartesian.
\end{sublem}

\textit{Proof:} We start by showing that the morphism
$$N(\mathcal{M}'(C,D)) \longrightarrow
N(\mathcal{M}'(C,D_{1}))\times^{h}_{N(\mathcal{M}'(C,D_{3}))}
N(\mathcal{M}'(C,D_{2}))$$
induces an injection on $\pi_{0}$ and an isomorphism
on all $\pi_{i}$ for $i>0$. For this, we consider the
induced diagram
of dg-categories
$$\xymatrix{
C\otimes D^{op} \ar[r]^-{u} \ar[d]_-{v} & C\otimes D^{op}_{1} \ar[d]^-{p} \\
C\otimes D^{op}_{2} \ar[r]_-{q} & C\otimes D^{op}_{3},}$$
where we keep the same names for the induced morphisms
after tensoring with $C$. It is then enough to show that
for $F$ and $G$ in $\mathcal{M}(C,D)$ the square
of path spaces
$$\xymatrix{
\Omega_{F,G}N(\mathcal{M}'(C,D)) \ar[r] \ar[d] &
\Omega_{u_{!}F,u_{!}G}N(\mathcal{M}'(C,D_{1})) \ar[d] \\
\Omega_{v_{!}F,v_{!}G}N(\mathcal{M}'(C,D_{2})) \ar[r] &
\Omega_{w_{!}F,w_{!}G}N(\mathcal{M}'(C,D_{3})),}$$
is homotopy cartesian (where $w=p\circ u$). Using
the natural equivalence between path spaces
in nerves of sub-categories of equivalences
in model categories
and mapping spaces of equivalences
(see \cite{dk}, and also \cite[Appendix A]{hagII}), one finds that
the previous diagram is in fact equivalent to the
following one
$$\xymatrix{
Map^{eq}(F,G) \ar[r] \ar[d] &
Map^{eq}(u_{!}F,u_{!}G) \ar[d] \\
Map^{eq}(v_{!}F,v_{!}G) \ar[r] &
Map^{eq}(w_{!}F,w_{!}G),}$$
where $Map^{eq}$ denotes the sub-simplicial set
of the mapping spaces consisting of all connected
components corresponding to equivalences.
By adjunction, this last diagram is equivalent to
$$\xymatrix{
Map^{eq}(F,G) \ar[r] \ar[d] &
Map^{eq}(F,u^{*}u_{!}G) \ar[d] \\
Map^{eq}(F,v^{*}v_{!}G) \ar[r] &
Map^{eq}(F,w^{*}w_{!}G).}$$
Therefore, to show that this last square is homotopy
cartesian, it is enough to prove that for any $G\in \mathcal{M}(C,D)$ the
natural morphism
$$G\longrightarrow u^{*}u_{!}G \times^{h}_{w^{*}w_{!}G}v^{*}v_{!}G$$
is an equivalence in $C \otimes D^{op}-Mod_{\mathbb{U}}$.
As this can be tested by fixing some object $x\in C$ and considering
the corresponding morphism
$$G(x,-)\longrightarrow (u^{*}u_{!}G \times^{h}_{w^{*}w_{!}G}v^{*}v_{!}G)(x,-)$$
in $D^{op}-Mod_{\mathbb{U}}$, we see that one can
assume that $C=\mathbf{1}$.
One can then write
$G=\underline{h}_{x}$ for some point $x\in D$. For $z\in D$, one has natural
isomorphisms
$$u^{*}u_{!}G(z)=D_{1}(u(z),u(x))
\qquad v^{*}v_{!}G(z)=D_{2}(v(z),v(x)) \qquad
w^{*}w_{!}G(z)=D_{3}(w(z),w(x)).$$
We therefore find that for any $z\in D$ the morphism
$$G(z)\longrightarrow (u^{*}u_{!}G \times^{h}_{w^{*}w_{!}G}v^{*}v_{!}G)(z)$$
can be written as
$$D(z,x)\longrightarrow D_{1}(u(z),u(x)) \times_{D_{3}(w(z),w(x))}^{h}
D_{2}(v(z),v(x)),$$
which by assumption on the morphism $p$ is a quasi-isomorphism
of complexes. This implies that
the morphism
$$G\longrightarrow u^{*}u_{!}G \times^{h}_{w^{*}w_{!}G}v^{*}v_{!}G$$
is an equivalence, and thus that
$$N(\mathcal{M}'(C,D)) \longrightarrow
N(\mathcal{M}'(C,D_{1}))\times^{h}_{N(\mathcal{M}'(C,D_{3}))}
N(\mathcal{M}'(C,D_{2}))$$
induces an injection on $\pi_{0}$ and an isomorphisms
on all $\pi_{i}$ for $i>0$. It only remains to show that the
above morphism is also surjective on connected components.

The set $\pi_{0}(N(\mathcal{M}'(C,D_{1}))\times^{h}_{N(\mathcal{M}'(C,D_{3}))}
N(\mathcal{M}'(C,D_{2})))$ can be described in the
following way. We consider
a category $\mathcal{N}$ whose objects are 5-tuples
$(F_{1},F_{2},F_{3};a,b)$, with
$F_{i}\in \mathcal{M}(C,D_{i})$ and where $a$ and $b$ are two morphisms
in $\mathcal{M}(C,D_{3})$
$$a : p_{!}(F_{1}) \longrightarrow F_{3} \longleftarrow
q_{!}(F_{2}) : b.$$
Morphisms in $\mathcal{N}$ are defined in the obvious way,
as morphisms $F_{i} \rightarrow G_{i}$ in
$\mathcal{M}(C,D_{i})$, commuting with the morphisms
$a$ and $b$. It is not hard
to check that $\pi_{0}(N(\mathcal{N}))$ is naturally
isomorphic to
$\pi_{0}(N(\mathcal{M}'(C,D_{1}))\times^{h}_{N(\mathcal{M}'(C,D_{3}))}
N(\mathcal{M}'(C,D_{2})))$. Furthermore, the natural
map
$$\pi_{0}(N(\mathcal{M}(C,D)))
\longrightarrow
\pi_{0}(N(\mathcal{N}))$$
is induced by the functor
$\mathcal{M}(C,D)
\longrightarrow \mathcal{N}$ that
sends an object
$F\in \mathcal{M}(C,D)$ to
$(u_{!}F,v_{!}F,w_{!}F;a,b)$ where
$a$ and $b$ are the two natural
isomorphisms
$$p_{!}u_{!}(F) \simeq w_{!}(F)\simeq q_{!}v_{!}(F).$$
Now, let $(F_{1},F_{2},F_{3};a,b) \in \mathcal{N}$,
and let us define an object $F\in Ho((C\otimes D^{op})-Mod_{\mathbb{U}})$
by the following formula
$$F:=u^{*}(F_{1})\times^{h}_{w^{*}(F_{3})}v^{*}(F_{1}).$$
Clearly, one has natural morphisms
in $Ho(C\otimes D_{i}^{op}-Mod_{\mathbb{U}})$
$$\mathbb{L}u_{!}(F)\rightarrow F_{1} \qquad \mathbb{L}v_{!}(F)\rightarrow F_{2}
\qquad \mathbb{L}w_{!}(F)\rightarrow F_{3}.$$
We claim that $F$ is right quasi-representable and that
these morphisms are in fact
isomorphisms. This will clearly finish the proof of
the surjectivity on connected components.
For this one can clearly assume that $C=\mathbf{1}$. One can then
write
$F_{i}=\underline{h}_{x_{i}}$, for some
$x_{i}\in D_{i}$. As $p$ is a fibration, the equivalence
$$a : p_{!}(\underline{h}_{x_{1}})=\underline{h}_{p(x_{1})} \longrightarrow
\underline{h}_{x_{3}}$$
can be lifted to an equivalence
$\underline{h}_{x_{1}} \longrightarrow \underline{h}_{x_{1}'}$
in $D^{op}_{1}-Mod$.
Replacing $x_{1}$ by $x_{1}'$ one can suppose that
$p(x_{1})=x_{3}$ and $a=id$. In the same way, the equivalence
$$b : q_{!}(\underline{h}_{x_{2}}) \longrightarrow\underline{h}_{p(x_{1})}$$
can be lifted to an equivalence
$\underline{h}_{x_{1}''} \longrightarrow \underline{h}_{x_{1}}$
in $D^{op}_{1}-Mod$.
Thus, replacing $x_{1}$ by $x_{1}''$ one can suppose that
$q(x_{2})=p(x_{1})=x_{3}$ and that $a$ and $b$ are the
identity morphisms. Then, clearly
$F\simeq \underline{h}_{x}$, where $x\in D$
is the point given by $(x_{1},x_{2},x_{3})$.
This shows that $F$ is right quasi-representable, and also that
the natural morphisms
$$u_{!}(F)\rightarrow F_{1} \qquad
v_{!}(F)\rightarrow F_{2}
 \qquad w_{!}(F)\rightarrow F_{3}$$
are equivalences.  \hfill $\Box$ \\

We have now finished the proof of sub-lemma \ref{sl5}
and thus of theorem \ref{t1}.  \hfill $\Box$ \\

Recall that $\mathcal{M}(Q(C),D)$ has been defined as the category
of equivalences between right quasi-representable
$Q(C)\otimes D^{op}$-modules $F$ such that $F(x,-)$ is cofibrant
in $D^{op}-Mod$ for any $x\in C$. This last condition is only
technical and useful for functorial reasons and does not affect the
nerve. Indeed, let $\mathcal{F}(Q(C),D)$ be the category of
all equivalences between right quasi-representable $(Q(C)\otimes D^{op})$-modules.
The natural inclusion functor
$$\mathcal{M}(Q(C),D) \longrightarrow \mathcal{F}(Q(C),D)$$
induces a weak equivalence on the corresponding nerves
as there exists a functor in the other direction just by taking
a cofibrant replacement (note that a cofibrant
$(Q(C)\otimes D^{op})$-module $F$ is such that
$F(x,-)$ is cofibrant for any $x\in Q(C)$, because of
Prop. \ref{p2}). In particular, theorem \ref{t1} implies the
existence of a string of weak equivalences
$$\xymatrix{
Map^{l}(C,D) \ar[r] & d(N(\mathcal{M}(\Gamma^{*}(C),D))) & \ar[l] N(\mathcal{M}(Q(C),D)) \ar[r] &
N(\mathcal{F}(Q(C),D)).}$$

The following corollary is a direct consequence of theorem
\ref{t1} and the above remark.

\begin{cor}\label{ct1}
Let $C$ and $D$ be two $\mathbb{U}$-small dg-categories.
Then, there exists a functorial
bijection between the set of maps $[C,D]$ in $Ho(dg-Cat_{\mathbb{U}})$, and
the set of isomorphism classes of right quasi-representable
objects in $Ho((C\otimes^{\mathbb{L}} D^{op})-Mod_{\mathbb{U}})$.
\end{cor}

Another important corollary of Theorem \ref{t1} is the
following.

\begin{cor}\label{ct1'}
Let $C$ be a $\mathbb{U}$-small dg-categories. Then, there exists
a functorial isomorphism between the set
$[\mathbf{1},C]$  and
the set of isomorphism classes of the category
$[C]$.
\end{cor}

\textit{Proof:} The Yoneda embedding
$\underline{h} : C \longrightarrow Int(C^{op}-Mod_{\mathbb{U}})$
induces a fully faithful functor
$$[C] \longrightarrow [Int(C^{op}-Mod_{\mathbb{U}})].$$
The essential image of this functor clearly is the
sub-category of quasi-representable $C^{op}$-modules. Therefore,
$[\underline{h}]$ induces a natural bijection between the
isomorphism classes of $[C]$ and the isomorphism classes
of quasi-representable objects in $[Int(C^{op}-Mod_{\mathbb{U}})]$.
As one has a natural equivalence
$[Int(C^{op}-Mod_{\mathbb{U}})]\simeq Ho(C^{op}-Mod_{\mathbb{U}})$
corollary \ref{ct1} implies the result.  \hfill $\Box$ \\

More generally, one can describe the higher homotopy groups
of the mapping spaces by the following formula.

\begin{cor}\label{cint}
Let $C$ be a $\mathbb{U}$-small dg-category, and $x\in C$ be an object.
Then, one has natural isomorphisms of groups
$$\pi_{1}(Map(\mathbf{1},C),x)\simeq Aut_{[C]}(x) \qquad
\pi_{i}(Map(\mathbf{1},C),x)\simeq H^{1-i}(C(x,x)) \; \forall \; i>1.$$
\end{cor}

\textit{Proof:} We use the general formula
$$\pi_{1}(N(W),x)\simeq Aut_{Ho(M)}(x) \qquad
\pi_{i}(N(W),x)\simeq \pi_{i-1}(Map_{M}(x,x),Id) \; \forall \; i>1,$$
for a model category $M$, its sub-category of equivalences $W$ and
a point $x\in M$ (see e.g. \cite[Cor. A.0.4]{hagII}). Applied to
$M=C^{op}-Mod_{\mathbb{U}}$ and using theorem \ref{t1} one finds
$$\pi_{1}(Map(\mathbf{1},C),x)\simeq Aut_{Ho(C^{op}-Mod)}(\underline{h}_{x}) \qquad
\pi_{i}(Map(\mathbf{1},C),x)\simeq \pi_{i-1}(Map_{C^{op}-Mod}(\underline{h}_{x},\underline{h}_{x}),Id) \; \forall \; i>1.$$
Using that the morphism $\underline{h}$ is quasi-fully faithful one finds
$$Aut_{Ho(C^{op}-Mod)}(\underline{h}_{x}) \simeq Aut_{[C]}(x) \qquad
\pi_{i-1}(Map_{C^{op}-Mod}(\underline{h}_{x},\underline{h}_{x}),Id)\simeq
H^{1-i}(C(x,x)).$$
\hfill $\Box$ \\

\begin{cor}\label{ct1''}
Let $C$ and $D$ be two $\mathbb{U}$-small dg-categories.
Let $Int((C\otimes^{\mathbb{L}}D^{op})-Mod_{\mathbb{U}}^{rqr})$
be the full sub-dg-category of
$Int((C\otimes^{\mathbb{L}}D^{op})-Mod_{\mathbb{U}})$ consisting
of all right quasi-representable objects. Then,
$Int((C\otimes^{\mathbb{L}}D^{op})-Mod_{\mathbb{U}}^{rqr})$
is isomorphic in $Ho(dg-Cat_{\mathbb{V}})$ to a $\mathbb{U}$-small
dg-category.
\end{cor}

\textit{Proof:} Indeed, we know by corollary \ref{ct1} that
the set of isomorphism classes of
$[Int((C\otimes^{\mathbb{L}}D^{op})-Mod_{\mathbb{U}}^{rqr})]$
is essentially $\mathbb{U}$-small, as it is
in bijection with $[C,D]$. Let us choose an essentially $\mathbb{U}$-small
full sub-dg-category $E$ in $Int((C\otimes^{\mathbb{L}}D^{op})-Mod_{\mathbb{U}}^{rqr})$ which
contains a set of representatives of isomorphism
classes of objects. As we already know that
the complexes of morphisms in $Int((C\otimes^{\mathbb{L}}D^{op})-Mod_{\mathbb{U}}^{rqr})$
are $\mathbb{U}$-small, the dg-category $E$ is
essentially $\mathbb{U}$-small, and thus isomorphic to a
$\mathbb{U}$-small dg-category. As $E$ is quasi-equivalent to
$Int((C\otimes^{\mathbb{L}}D^{op})-Mod_{\mathbb{U}}^{rqr})$
this implies the result.
\hfill $\Box$ \\

We finish by the following last corollary.

\begin{cor}\label{cint2}
Let $C$ and $D$ be two $\mathbb{U}$-small dg-categories, and
let $f,g : C \longrightarrow D$ be two morphisms
with corresponding $(C\otimes^{\mathbb{L}} D^{op})$-modules
$\phi(f)$ and $\phi(g)$. Then, there exists a natural weak equivalence
of simplicial sets
$$\Omega_{f,g}Map_{dg-Cat}(C,D)\simeq Map_{(C\otimes^{\mathbb{L}} D^{op})-Mod}^{eq}(\phi(f),\phi(g)),$$
where $Map^{eq}(\phi(f),\phi(g))$ is the sub-simplicial set of
$Map(\phi(f),\phi(g))$ consisting of equivalences.
\end{cor}

\textit{Proof:} This follows immediately from
theorem \ref{t1} and the standard relations between
path spaces of nerves of equivalences in a model category and
its mapping spaces (see e.g. \cite[Appendix A]{hagII}). \hfill $\Box$ \\

\section{The simplicial structure}

Let $K\in SSet_{\mathbb{U}}$ be a $\mathbb{U}$-small
simplicial set and $C\in dg-Cat_{\mathbb{U}}$. One can form
the derived tensor product $K\otimes^{\mathbb{L}}C \in Ho(dg-Cat_{\mathbb{U}})$,
as well as the derived exponential $C^{\mathbb{R}K}$. One has the usual
adjunction isomorphism
$$[K\otimes^{\mathbb{L}}C,D]\simeq [C,D^{\mathbb{R}K}]\simeq [K,Map(C,D)].$$
Let $\Delta(K)$ be the
simplex category  of $K$. An object of $\Delta(K)$
is therefore a pair $(n,a)$ with $n\in \Delta$ and $x\in K_{n}$.
A morphism $(n,x) \rightarrow (m,y)$ is the data
of a morphism $u : [n] \rightarrow [m]$ in $\Delta$
such that $u^{*}(y)=x$. The simplicial set
$K$ is then naturally weakly equivalent to the
homotopy colimit of the constant diagram
$$\Delta(K) \longrightarrow * \in SSet.$$
In other words, one has a natural weak equivalence
$$N(\Delta(K))\simeq K.$$
We now consider
$\Delta(K)_{k}$ the $k$-linear category
freely generated by the category $\Delta(K)$, and
consider $\Delta(K)_{k}$ as an object
in $dg-Cat_{\mathbb{U}}$.

\begin{thm}\label{p5}
Let $C$ and $D$ be two $\mathbb{U}$-small
dg-categories, and $K\in SSet_{\mathbb{U}}$.
Then, there exists a functorial injective map
$$[K\otimes^{\mathbb{L}}C,D] \longrightarrow [\Delta(K)_{k}\otimes^{\mathbb{L}}C,D].$$
Moreover, the image of this map consists exactly of all
morphism $\Delta(K)_{k}\otimes^{\mathbb{L}}C\longrightarrow D$
in $Ho(dg-Cat_{\mathbb{U}})$ such that for any
$c\in C$ the induced functor
$$\Delta(K)_{k} \longrightarrow [D]$$
sends all morphisms in $\Delta(K)_{k}$ to
isomorphisms in $[D]$.
\end{thm}

\textit{Proof:} Using our theorem \ref{t1} one finds natural
equivalences
$$[K\otimes^{\mathbb{L}}C,D]\simeq [K,Map(C,D)]\simeq
[K,N(\mathcal{M}(Q(C),D))].$$

We then use the next technical lemma.

\begin{lem}\label{lp5}
Let $M$ be a $\mathbb{V}$-small $\mathbb{U}$-combinatorial model category
and $K\in SSet_{\mathbb{U}}$. Let $W\subset M$ be the sub-category
of equivalences in $M$. Then, there exists a natural bijection
between $[K,N(W)]_{SSet_{\mathbb{V}}}$ and the set of isomorphism
classes of objects
$F \in Ho(M^{\Delta(K)})$ corresponding to functors $F : \Delta(K)
\longrightarrow M$
sending all morphisms of $\Delta(K)$ to
equivalences in $M$.
\end{lem}

\textit{Proof:} First of all, the lemma is invariant by
changing $M$ up to a Quillen equivalence, and thus
by \cite{du} one can suppose that
$M$ is a simplicial model category.
The proof of the lemma will use some techniques of
simplicial localizations \`a la Dwyer-Kan,
as well as some result about $S$-categories. We start by
a short digression on the subject.

We recall the existence of a model category of
$S$-categories, as shown in \cite{be}, and which is similar
to the one we use on dg-categories. This model category
will be denoted by $S-Cat$ (or $S-Cat_{\mathbb{V}}$ if one needs
to specify the universe). For any
$\mathbb{V}$-small category $C$ with a sub-category $S\subset C$, one can
form a $\mathbb{V}$-small $S$-category $L(C,S)$ by formally inverting the
morphisms in $S$ in a homotopy meaningful way (see e.g. \cite{dk2}).
Using the language of model categories, this means that for any
$\mathbb{V}$-small $S$-category $T$, there exists functorial isomorphisms
between $[L(C,S),T]_{S-Cat}$ and the subset of $[C,T]_{S-Cat}$ consisting
of all morphisms sending $S$ to
isomorphisms in $[T]$ (the category $[T]$ is defined by taking connected
component of simplicial sets of morphisms in $T$). Finally, one can define
a functor $N : Ho(S-Cat_{\mathbb{V}}) \longrightarrow Ho(SSet_{\mathbb{V}})$
by sending an $S$-category to its nerve. It is well known that
the functor $N$ becomes an equivalence when restricted to
$S$-categories $T$ such that $[T]$ is a groupoid (this is
just another way to state delooping theory). Finally, for any
category $C$ with a sub-category $S\subset C$, one has a natural
weak equivalence $N(L(C,S))\simeq N(C)$.

Now, as explained in \cite[Prop. A.0.6]{hagII},
$N(W)$ can be also interpreted as the nerve
of the $S$-category $\mathcal{G}(M)$, of cofibrant and
fibrant objects in $M$ together with their simplicial sets of equivalences.
One therefore has natural isomorphism
$$[K,N(W)]\simeq [N(\Delta(K)),N(\mathcal{G}(M))]\simeq [L(\Delta(K),\Delta(K)),\mathcal{G}(M)].$$
Furthermore, as all morphisms in $[\mathcal{G}(M)]$ are isomorphisms one finds
a bijection between $[K,N(W)]$ and $[\Delta(K),\mathcal{G}(M)]$.
Let $Int(M)$ be the $S$-category of fibrant and cofibrant objects
in $M$ together with their simplicial sets of morphisms.
Then, as $\mathcal{G}(M)$ is precisely the sub-$S$-category of
$Int(M)$ consisting of equivalences, the set $[\Delta(K),\mathcal{G}(M)]$ is also
the subset of $[\Delta(K),Int(M)]$ consisting of all morphisms such that the induced
functor $\Delta(K) \longrightarrow [Int(M)]\simeq Ho(M)$
sends all morphisms to isomorphisms.
Finally,
it turns out that the same results as
our lemmas \ref{sl2} and \ref{sl1} are valid in the context of $S$-categories (their proofs
are exactly the same). Therefore, we see that $[\Delta(K),Int(M)]$ is in a natural bijection with
isomorphism classes of objects in $Ho(M^{\Delta(K)})$. Putting all of this together
gives the lemma. \hfill $\Box$ \\

We apply the previous lemma to the case where
$M:=(C\otimes^{\mathbb{L}}D^{op})-Mod_{\mathbb{U}}$, and we find a
natural injection $[K,N(W)] \hookrightarrow
Iso(Ho(M^{\Delta(K)}))$, whose image consists of all functors
$\Delta(K) \rightarrow M$ sending all morphisms of $\Delta(K)$ to
equivalences in $M$. Composing with the natural inclusion
$\mathcal{M}(Q(C),D)\subset M$ provides a natural injection of
$$[K,N(\mathcal{M}(Q(C),D))] \subset [K,N(W)] \subset
Iso(Ho(M^{\Delta(K)})).$$
By the construction of the bijection of lemma \ref{lp5} one easily sees
that the image of this inclusion consists of all
functors $F : \Delta(K) \longrightarrow W$ such that
for any $k\in K$ one has $F(k)\in \mathcal{M}(Q(C),D)$.
Finally, one clearly has a natural equivalence of categories, compatible
with the model structures
$$M^{\Delta(K)}\simeq (C\otimes^{\mathbb{L}}D^{op})-Mod_{\mathbb{U}}^{\Delta(K)_{k}}\simeq
(\Delta(K)_{k}\otimes C\otimes^{\mathbb{L}}D^{op})-Mod_{\mathbb{U}},$$
inducing a bijection between
$Iso(Ho(M^{\Delta(K)}))$ and the isomorphism classes of
objects in $Ho((\Delta(K)_{k}\otimes C\otimes^{\mathbb{L}}D^{op})-Mod_{\mathbb{U}})$.
Another application of theorem \ref{t1} easily implies the result.
\hfill $\Box$ \\

\section{Internal Hom's}

Let us recall that $Ho(dg-Cat_{\mathbb{U}})$ is endowed with
the symmetric monoidal structure $\otimes^{\mathbb{L}}$.
Recall that the monoidal structure $\otimes^{\mathbb{L}}$
is said to be closed if for any two objects
$C$ and $D$ in $Ho(dg-Cat_{\mathbb{U}})$ the functor
$A \mapsto [A\otimes^{\mathbb{L}}C,D]$ is representable
by an object $\mathbb{R}\underline{Hom}(C,D)\in Ho(dg-Cat_{\mathbb{U}})$.
Recall also from
corollary \ref{ct1''} that the $\mathbb{V}$-small dg-category
$Int((C\otimes^{\mathbb{L}}D^{op})-Mod_{\mathbb{U}}^{rqr})$
is essentially $\mathbb{U}$-small and therefore can be
considered as an object in $Ho(dg-Cat_{\mathbb{U}})$.

\begin{thm}\label{t2}
The monoidal category $(Ho(dg-Cat_{\mathbb{U}}),\otimes^{\mathbb{L}})$
is closed. Furthermore, for any two $\mathbb{U}$-small dg-categories
$C$ and $D$ one has a natural isomorphism in $Ho(dg-Cat_{\mathbb{U}})$
$$\mathbb{R}\underline{Hom}(C,D)\simeq
Int((C\otimes^{\mathbb{L}}D^{op})-Mod_{\mathbb{U}}^{rqr}).$$
\end{thm}

\textit{Proof:} The proof is essentially the same
as for theorem \ref{t1} and is also based on the same
lemmas \ref{sl2} and \ref{sl1}. Indeed, from these
two lemmas one extracts the following result.

\begin{lem}\label{l3}
Let $M$ be $C(k)_{\mathbb{U}}$-enriched
$\mathbb{U}$-cofibrantly generated model category
which is $\mathbb{V}$-small. We assume that
the domain and codomain of a set of generating
cofibrations are cofibrant in $M$.
 Let $M_{0}$ be
a full sub-category of $M$ which is closed by equivalences, and
$Int(M_{0})$ be the full sub-dg-category of
$Int(M)$ consisting of all objects belonging to $M_{0}$.
Let $A$ be a cofibrant and
$\mathbb{U}$-small dg-category, and let $Ho(M_{0}^{A})$ be the
full sub-category of $Ho(M^{A})$ consisting of objects
$F\in Ho(M^{A})$ such that $F(a)\in M_{0}$ for any $a\in A$.
Then, one has a natural
isomorphism
$$\phi : [A,Int(M_{0})]\simeq Iso(Ho(M_{0}^{A})).$$
\end{lem}

\textit{Proof:} The morphism
$$\phi : [A,Int(M_{0})]\longrightarrow Iso(Ho(M_{0}^{A}))$$
simply sends a morphism $A \longrightarrow
Int(M_{0})$ to the corresponding object in $M_{0}^{A}$.
Using our proposition \ref{p1} it is easy to see that
this maps sends homotopic morphisms to isomorphic
objects in $Ho(M_{0}^{A})$, and is therefore well defined.
As for the proof of lemma \ref{sl2}, the morphism
$\phi$ is clearly surjective. Let $u,v : A \longrightarrow
Int(M_{0})$ be two morphisms of dg-categories such that
the corresponding objects in $Ho(M_{0}^{A})$ are isomorphic.
Then, these objects are isomorphic in $Ho(M^{A})$, which
implies by lemma \ref{sl1} that the two compositions
$$u',v' : A \longrightarrow Int(M_{0}) \longrightarrow Int(M)$$
are homotopic in $dg-Cat_{\mathbb{V}}$. Let
$$\xymatrix{
A \ar[rd]^-{u'} \ar[d] & \\
A' \ar[r]^-{H} & Int(M) \\
A \ar[ru]_-{v'} \ar[u]}$$
be a homotopy between $u'$ and $v'$. As $M_{0}$ is closed by equivalences
in $M$ one clearly sees that the morphism $H$ factors through the
sub-dg-category $Int(M_{0})$, showing that
$u$ and $v$ are homotopic. \hfill $\Box$ \\

We come back to the proof of theorem \ref{t2}.
Using our theorem \ref{t1} one has a natural
isomorphism
$$[A\otimes^{\mathbb{L}}C,D]\simeq
Iso(Ho(((A\otimes^{\mathbb{L}} C)\otimes^{\mathbb{L}} D^{op})-Mod_{\mathbb{U}}^{rqr}))\simeq
Iso(Ho(((C\otimes^{\mathbb{L}} D^{op})-Mod_{\mathbb{U}}^{rqr})^{A})).$$
An application of lemma \ref{l3} (with
$M=(C\otimes^{\mathbb{L}}D^{op})-Mod_{\mathbb{U}}$ and
$M_{0}$ the full sub-category of right quasi-representable
objects)
shows that
one has a natural isomorphism
$$[A,Int((C\otimes^{\mathbb{L}}D^{op})-Mod_{\mathbb{U}}^{rqr})]\simeq
Iso(Ho(((C\otimes^{\mathbb{L}} D^{op})-Mod_{\mathbb{U}}^{rqr})^{A})).$$
Putting this together one finds a natural isomorphism
$$[A\otimes^{\mathbb{L}}C,D]\simeq [A,Int((C\otimes^{\mathbb{L}}D^{op})-Mod_{\mathbb{U}}^{rqr})]$$
showing the theorem. \hfill $\Box$ \\

\begin{cor}\label{cp5}
For any $C$ and $D$ two $\mathbb{U}$-small
dg-categories, and any $K\in SSet_{\mathbb{U}}$, one has
a functorial isomorphism in $Ho(dg-Cat_{\mathbb{U}})$
$$K\otimes^{\mathbb{L}}(C\otimes^{\mathbb{L}}D)\simeq
(K\otimes^{\mathbb{L}}C)\otimes^{\mathbb{L}}D.$$
\end{cor}

\textit{Proof:} This follows easily from Thm. \ref{p5},
Thm. \ref{t2} and the Yoneda lemma applied to $Ho(dg-Cat_{\mathbb{U}})$.
\hfill $\Box$ \\

\begin{cor}\label{cp5'}
For any $C$, $D$ and $E$ three $\mathbb{U}$-small
dg-categories one has
a functorial isomorphism in $Ho(SSet_{\mathbb{U}})$
$$Map(C\otimes^{\mathbb{L}}D,E)\simeq
Map(C,\mathbb{R}\underline{Hom}(D,E)).$$
\end{cor}

\textit{Proof:} By Cor. \ref{cp5}, for any $K\in SSet_{\mathbb{U}}$, one has
functorial isomorphisms
$$[K,Map(C\otimes^{\mathbb{L}}D,E)]\simeq
[K\otimes^{\mathbb{L}}(C\otimes^{\mathbb{L}}D),E]\simeq
[(K\otimes^{\mathbb{L}}C)\otimes^{\mathbb{L}}D,E]\simeq$$
$$[K\otimes^{\mathbb{L}}C,\mathbb{R}\underline{Hom}(D,E)]\simeq
[K,Map(C,\mathbb{R}\underline{Hom}(D,E))].$$
\hfill $\Box$ \\

\begin{cor}\label{cp5''}
Let $C\in dg-Cat_{\mathbb{U}}$ be
a dg-category. Then the functor
$$-\otimes^{\mathbb{L}}C : dg-Cat_{\mathbb{U}} \longrightarrow dg-Cat_{\mathbb{U}}$$
commutes with homotopy colimits.
\end{cor}

\textit{Proof:} This follows formally from Cor. \ref{cp5'}. \hfill $\Box$ \\

\begin{cor}\label{cp5'''}
Let $C \longrightarrow D$ be a quasi-fully faithful morphism in $dg-Cat_{\mathbb{U}}$. Then, for
any $B\in dg-Cat_{\mathbb{U}}$ the induced morphism
$$\mathbb{R}\underline{Hom}(B,C) \longrightarrow \mathbb{R}\underline{Hom}(B,D)$$
is quasi-fully faithful.
\end{cor}

\textit{Proof:} Using Lem. \ref{lmono} it is enough to show that
$\mathbb{R}\underline{Hom}(B,-)$ preserves homotopy monomorphisms. But this
follows formally from Cor. \ref{cp5'}. \hfill $\Box$ \\

\section{Morita morphisms and bi-modules}

In this paragraph we will use the following notations. For any
$C\in dg-Cat_{\mathbb{U}}$ one sets
$$\widehat{C}:=Int(C^{op}-Mod_{\mathbb{U}})\in dg-Cat_{\mathbb{V}}.$$
By theorem \ref{t2} and lemma \ref{l3}, one has an isomorphism in
$Ho(dg-Cat_{\mathbb{V}})$
$$\widehat{C}\simeq \mathbb{R}\underline{Hom}(C^{op},Int(C(k)_{\mathbb{U}})) \in Ho(dg-Cat_{\mathbb{V}}).$$
Indeed, lemma \ref{l3} implies that for any $A\in dg-Cat_{\mathbb{U}}$ one has
$$[A,\widehat{C}]\simeq Iso(Ho((A\otimes^{\mathbb{L}}C^{op})-Mod_{\mathbb{U}}))\simeq
[A\otimes^{\mathbb{L}}C^{op},\widehat{\mathbf{1}}].$$
Note also that
$$Int(C(k)_{\mathbb{U}})\simeq \widehat{\mathbf{1}}.$$
We will also consider $\widehat{C}_{pe}$ the full sub-dg-category of
$\widehat{C}$ consisting of $C^{op}$-modules which are
homotopically finitely presented. In other words,
a $C^{op}$-module $F$ is in $\widehat{C}_{pe}$ if for any
filtered diagram of objects $G_{i}$ in $C^{op}-Mod_{\mathbb{U}}$, the natural morphism
$$Colim_{i}Map(F,G_{i}) \longrightarrow Map(F,Colim_{i}G_{i})$$
is a weak equivalence. It is easy to check that the
objects in $\widehat{C}_{pe}$ are precisely the objects equivalent to retracts of finite
cell $C^{op}$-modules. To be more precise, an object $F\in Ho(\widehat{C})$ is
in $Ho(\widehat{C}_{pe})$ if and only if it is a retract in $Ho(\widehat{C})$ of
an object $G$ for which there exists a finite sequence of morphisms of $C^{op}$-modules
$$\xymatrix{
0 \ar[r] & G_{1} \ar[r] & G_{2} \ar[r] & \dots \ar[r] & G_{n}=G,}$$
in such a way that for any $i$ there exists a push-out square
$$\xymatrix{G_{i} \ar[r] & G_{i+1} \\
A\otimes \underline{h}_{x} \ar[u] \ar[r] & \ar[u] B\otimes \underline{h}_{x}}$$
for some $x\in C$, and some cofibration $A\rightarrow B$ in $C(k)$ with
$A$ and $B$ bounded complexes of projective modules of finite type.

Objects in $\widehat{C}_{pe}$ will also be called
\emph{compact} or \emph{perfect} (note that they are precisely
the compact objects in the triangulated category $[\widehat{C}]$, in the
usual sense). More generally, for any dg-category $T$, we will write
$T_{pe}$ for the full sub-dg-category of $T$ consisting of compact objects
(i.e. the objects $x$ such that $[T](x,-)$ commutes with
(infinite) direct sums). \\

Let us consider $C$ and $D$ two $\mathbb{U}$-small dg-categories,
and $u : \widehat{C} \longrightarrow \widehat{D}$ a morphism in $Ho(dg-Cat_{\mathbb{V}})$.
Then, $u$ induces a functor, well defined up to an (non-unique) isomorphism
$$[u] : [\widehat{C}]\longrightarrow [\widehat{D}].$$
We will say that the morphism $u$ is \emph{continuous} if the functor
$[u]$ commutes with $\mathbb{U}$-small direct sums. Note that
$[\widehat{C}]$ and $[\widehat{D}]$ are the homotopy categories of the model
categories of $C^{op}$-modules and $D^{op}$-modules, and thus these two
categories always have direct sums. More generally, we will denote by
$\mathbb{R}\underline{Hom}_{c}(\widehat{C},\widehat{D})$ the full sub-dg-category
of $\mathbb{R}\underline{Hom}(\widehat{C},\widehat{D})$ consisting of
continuous morphisms.

\begin{df}\label{d5}
Let $C$ and $D$ be two $\mathbb{U}$-small dg-categories.
\begin{enumerate}
\item
The dg-category of \emph{Morita morphisms} from $C$ to $D$ is
$\mathbb{R}\underline{Hom}_{c}(\widehat{C},\widehat{D})$.
\item The dg-category of \emph{perfect Morita morphisms} from $C$ to $D$ is
$\mathbb{R}\underline{Hom}(\widehat{C}_{pe},\widehat{D}_{pe})$.
\end{enumerate}
\end{df}

We warn the reader that there are in general no relations between
the dg-category $\mathbb{R}\underline{Hom}(\widehat{C}_{pe},\widehat{D}_{pe})$ and
$\mathbb{R}\underline{Hom}_{c}(\widehat{C},\widehat{D})_{pe}$. An example where these
two objects agree will be given in Thm. \ref{tfour2}. \\

\begin{thm}\label{t3}
Let $C\in dg-Cat_{\mathbb{U}}$, and let us consider the Yoneda
embedding $\underline{h} : C \longrightarrow \widehat{C}$. Let
$D$ be any $\mathbb{U}$-small dg-category.
\begin{enumerate}
\item
The pull-back functor
$$\underline{h}^{*} : \mathbb{R}\underline{Hom}_{c}(\widehat{C},\widehat{D}) \longrightarrow
\mathbb{R}\underline{Hom}(C,\widehat{D})$$
is an isomorphism in $Ho(dg-Cat_{\mathbb{V}})$.
\item The pull-back functor
$$\underline{h}^{*} : \mathbb{R}\underline{Hom}(\widehat{C}_{pe},\widehat{D}_{pe}) \longrightarrow
\mathbb{R}\underline{Hom}(C,\widehat{D}_{pe})$$
is an isomorphism in $Ho(dg-Cat_{\mathbb{V}})$.
\end{enumerate}
\end{thm}

\textit{Proof:} We start by proving $(1)$. \\

Using the universal properties of internal Hom's one reduces the problem to show that for
any $A\in dg-Cat_{\mathbb{U}}$, the morphism\footnote{We prefer to change
notation from $\underline{h}$ to $l$ during the proof, just
in order to avoid future confusions.}
$$l:=\underline{h} : C \longrightarrow \widehat{C}$$
induces a bijective morphism
$$l^{*} : [\widehat{C}\otimes^{\mathbb{L}} A,\widehat{D}]_{c} \longrightarrow
[C\otimes^{\mathbb{L}} A,\widehat{D}],$$
where by definition $[\widehat{C}\otimes^{\mathbb{L}} A,\widehat{D}]_{c}$
is the subset of $[\widehat{C}\otimes^{\mathbb{L}} A,\widehat{D}]$ consisting
of morphisms $f : \widehat{C}\otimes^{\mathbb{L}} A \longrightarrow \widehat{D}$
such that for any object $a\in A$ the induced morphism $f(-,a) : \widehat{C} \longrightarrow 
\widehat{D}$ is continuous. Now, as $\widehat{D}=\mathbb{R}\underline{Hom}(D^{op},\widehat{\mathbf{1}})$, 
one has natural bijections
$$[C\otimes^{\mathbb{L}} A,\widehat{D}]\simeq [C,\widehat{A^{op}\otimes^{\mathbb{L}}D}] \qquad
[\widehat{C}\otimes^{\mathbb{L}} A,\widehat{D}]_{c} \simeq [\widehat{C},\widehat{A^{op}\otimes^{\mathbb{L}}D}]_{c}.$$ 
Therefore, we have to prove that for any $A$ the induced morphism
$$l^{*} : [\widehat{C},\widehat{A^{op}\otimes^{\mathbb{L}}D}]_{c} \longrightarrow
[C,\widehat{A^{op}\otimes^{\mathbb{L}}D}],$$
is bijective. For this, we consider the
quasi-fully faithful morphism in $dg-Cat_{\mathbb{W}}$
for some universe $\mathbb{V}\in \mathbb{W}$
$$\widehat{A^{op}\otimes^{\mathbb{L}}D}\simeq Int((A\otimes^{\mathbb{L}}D^{op})-Mod_{\mathbb{U}}) \longrightarrow
\widehat{A^{op}\otimes^{\mathbb{L}}D}_{\mathbb{V}}:=
Int((A\otimes^{\mathbb{L}}D^{op})-Mod_{\mathbb{V}}).$$
One has a commutative square
$$\xymatrix{
[\widehat{C},\widehat{A^{op}\otimes^{\mathbb{L}}D}]_{c} \ar[r] \ar[d] &
[\widehat{C},\widehat{A^{op}\otimes^{\mathbb{L}}D}_{\mathbb{V}}]_{c} \ar[d] \\
[C,\widehat{A^{op}\otimes^{\mathbb{L}}D}] \ar[r] &
[C,\widehat{A^{op}\otimes^{\mathbb{L}}D}_{\mathbb{V}}].}$$

We claim that the right vertical morphism is bijective. For this,
we use lemma \ref{l3} which implies that it is enough to show
the following lemma.

\begin{lem}\label{lt3}
Let $C$ be a $\mathbb{U}$-small dg-category and
$M$ a $\mathbb{V}$-combinatorial $C(k)_{\mathbb{V}}$-model category
which is $\mathbb{W}$-small for some $\mathbb{V}\in \mathbb{W}$.
We assume that the domain and codomain of a set of
generating cofibrations are cofibrant in $M$.
We also assume that for any cofibrant object $X\in M$, and any
quasi-isomorphism $Z\longrightarrow Z'$ in $C(k)$, the
induced morphism
$$Z\otimes X \longrightarrow Z'\otimes X$$
is an equivalence in $M$.
Then, the Quillen adjunction
$$l_{!} : M^{C} \longrightarrow M^{\widehat{C}} \qquad
M^{C} \longleftarrow M^{\widehat{C}} : l^{*}$$
induces a fully faithful functor
$$\mathbb{L}l_{!} : Ho(M^{C}) \longrightarrow
Ho(M^{\widehat{C}})$$
whose essential image consists of all $\widehat{C}$-modules
corresponding to continuous morphisms in $Ho(dg-Cat_{\mathbb{W}})$.
\end{lem}

\textit{Proof:} First of all, the modules $F \in Ho(M^{\widehat{C}})$
corresponding to continuous morphisms are precisely the ones
for which for any $\mathbb{U}$-small family of objects $x_{i}\in \widehat{C}$, the
natural morphism
$$\bigoplus^{\mathbb{L}} F(x_{i}) \longrightarrow F(\oplus_{i}x_{i})$$
is an isomorphism in $Ho(M)$.

We start by showing that $\mathbb{L}l_{!}$ is fully faithful.
As both functors $\mathbb{L}l_{!}$ and $l^{*}$ commute
with homotopy colimits, it is enough to show that
for any $x\in C$ and any  $X\in M$, the adjunction morphism
$$X\otimes^{\mathbb{L}} \underline{h}^{x} \longrightarrow l^{*}\mathbb{L}l_{!}(
X\otimes^{\mathbb{L}} \underline{h}^{x})$$
is an isomorphism in $Ho(M^{C})$. But this follows immediately from the fact that
the morphism of dg-categories $l$ is fully faithful
and our hypothesis on $M$.

It remains to show that for any $F\in Ho(M^{\widehat{C}})$, corresponding to
a continuous morphism, the adjunction morphism
$$\mathbb{L}l_{!}l^{*}(F) \longrightarrow F$$
is an isomorphism in $Ho(M^{\widehat{C}})$. As we already
know that $\mathbb{L}l_{!}$ is fully faithful it is enough to
show that the functor $l^{*}$ is conservative when restricted
to the sub-category of modules corresponding to continuous functors.
Let $u : F \longrightarrow G$ be morphism between such modules, and let us
assume that $l^{*}(F) \longrightarrow l^{*}(G)$
is an isomorphism in $Ho(M^{C})$. We need to show that
$u$ itself is an isomorphism in $Ho(M^{\widehat{C}})$.

\begin{sublem}\label{sllt3}
Let $F : \widehat{C} \longrightarrow M$ be a morphism of
dg-categories corresponding to a continuous morphism.
\begin{enumerate}
\item
Let $X : I \longrightarrow C^{op}-Mod_{\mathbb{U}}$
be a $\mathbb{U}$-small diagram of cofibrant objects in
$C^{op}-Mod_{\mathbb{U}}$. Then, the natural morphism
$$Hocolim_{i}F(X_{i}) \longrightarrow
F(Hocolim_{i}X_{i})$$
is an isomorphism in $Ho(M)$.
\item Let $Z\in C(k)_{\mathbb{U}}$ and $X\in M$. Then, the natural
morphism
$$Z\otimes^{\mathbb{L}}F(X) \longrightarrow F(Z\otimes^{\mathbb{L}}X)$$
is an isomorphism in $Ho(M)$.
\end{enumerate}
\end{sublem}

\textit{Proof of sub-lemma \ref{sllt3}:}
$(1)$ As any homotopy colimit is a composition of
homotopy push-outs and infinite (homotopy) sums, it is enough to check the
sub-lemma for one of these colimits. For the direct sum case this is
our hypothesis on $F$. It remains to show that $F$ commutes with
homotopy push-outs. For this we assume that $F$ is fibrant and cofibrant, and thus
is given by a morphism of dg-categories $\widehat{C} \longrightarrow Int(M)$.

We consider the commutative diagram of dg-categories
$$\xymatrix{
(\widehat{C})^{op} \ar[r]^-{F} \ar[d] & Int(M)^{op} \ar[d] \\
Int(\widehat{C}-Mod_{\mathbb{V}}) \ar[r]_-{F_{!}} &
Int(Int(M)-Mod_{\mathbb{V}}),}$$
where the vertical morphisms are the dual Yoneda embeddings $\underline{h}^{(-)}$.
The functor $F_{!}$ being left Quillen clearly commutes, up to equivalences,
with homotopy push-outs. Furthermore, as the model categories
$\widehat{C}-Mod_{\mathbb{V}}$ and $Int(M)-Mod_{\mathbb{V}}$ are stable model categories, this
implies that $F_{!}$ also commutes, up to equivalence, with homotopy pull-backs.
Furthermore, the morphism $\underline{h}^{(-)}$ sends
homotopy push-out squares to homotopy pull-back squares, and moreover
a square in $Int(M)$ is a homotopy push-out square if and only if
its image by $\underline{h}$ is a homotopy pull-back square
in $Int(M)-Mod_{\mathbb{V}}$. We deduce from these remarks that
$F$ preserves homotopy push-out squares.  \\

$(2)$ Any complex $Z$ can be constructed from the trivial
complex $k$ using homotopy colimits and loop objects.
As we already know that $F$ commutes with homotopy colimits, it
is enough to see that it also commutes with loop objects.
But the loop functor is inverse, up to equivalence, to
the suspension functor. The suspension being a homotopy push-out,
$F$ commutes with it, and therefore $F$ commutes with the loop functor.
\hfill $\Box$ \\

Now, let us come back to our morphism
$u : F \longrightarrow G$ such that $l^{*}(u)$ is
an equivalence. Let $X$ be an object in $\widehat{C}$. We know that
$X$ can be written as the homotopy colimit of
objects of the form $Z\otimes^{\mathbb{L}}\underline{h}_{x}$
with $x\in C$ and $Z\in C(k)$.
Therefore, one has a commutative
diagram in $Ho(M)$
$$\xymatrix{
Hocolim_{i} F(Z_{i}\otimes^{\mathbb{L}}\underline{h}_{x_{i}}) \ar[d] \ar[r]^-{u}
& Hocolim_{i} G(Z_{i}\otimes^{\mathbb{L}}\underline{h}_{x_{i}}) \ar[d] \\
F(X) \ar[r]^-{u} & G(X).}$$
By the sub-lemma $(1)$ the vertical morphisms are isomorphisms in $Ho(M)$, and
the top horizontal morphism is also by hypothesis
and the sub-lemma $(2)$.
Thus, the bottom horizontal morphism is an isomorphism in $Ho(M)$, and this
for any $X\in \widehat{C}$. This shows that $l^{*}$ is conservative when
restricted to continuous morphisms, and thus finishes the proof of the
lemma  \ref{lt3}. \hfill $\Box$ \\

We come back to our commutative diagram
$$\xymatrix{
[\widehat{C},\widehat{A^{op}\otimes^{\mathbb{L}}D}]_{c} \ar[r] \ar[d] &
[\widehat{C},\widehat{A^{op}\otimes^{\mathbb{L}}D}_{\mathbb{V}}]_{c} \ar[d] \\
[C,\widehat{A^{op}\otimes^{\mathbb{L}}D}] \ar[r] &
[C,\widehat{A^{op}\otimes^{\mathbb{L}}D}_{\mathbb{V}}].}$$
Lemma \ref{lt3} shows that the right vertical
morphism is bijective, and corollary \ref{clmono} implies that
the horizontal morphisms are injective. It remains to show that
a morphism $u\in [\widehat{C},\widehat{A^{op}\otimes^{\mathbb{L}}D}_{\mathbb{V}}]_{c}$,
whose restriction $C \longrightarrow \widehat{A^{op}\otimes^{\mathbb{L}}D}_{\mathbb{V}}$
factors thought
$\widehat{A^{op}\otimes^{\mathbb{L}}D}$, itself factors through
$\widehat{A^{op}\otimes^{\mathbb{L}}D}$. But this is true
as by sub-lemma \ref{sllt3} the image by $u$ of any
$C^{op}$-module can be written as a $\mathbb{U}$-small homotopy colimit
of objects of the form $Z\otimes^{\mathbb{L}}u(l(x)) $
for $Z\in C(k)_{\mathbb{U}}$ and $x\in C$. Therefore,
if the restriction of $u$ to $C$ has $\mathbb{U}$-small images, then so does
$u$ itself. This finishes the proof of theorem \ref{t3} $(1)$. \\

$(2)$ We consider the quasi-fully faithful morphism
$\widehat{D}_{pe} \longrightarrow \widehat{D}$. We therefore have
a homotopy commutative diagram
$$\xymatrix{
\mathbb{R}\underline{Hom}(\widehat{C}_{pe},\widehat{D}_{pe}) \ar[r]
\ar[d] & \ar[d] \mathbb{R}\underline{Hom}(\widehat{C}_{pe},\widehat{D}) \\
\mathbb{R}\underline{Hom}(C,\widehat{D}_{pe})  \ar[r] &  \mathbb{R}\underline{Hom}(C,\widehat{D}),}$$
where the horizontal morphisms are quasi-fully faithful by Cor. \ref{cp5'''}. We claim that the
right vertical morphism is a quasi-equivalence. For this, using the universal
properties of internal Hom's, it is enough to show that the induced morphism
$$[\widehat{C}_{pe},\widehat{D}] \longrightarrow [C,\widehat{D}]$$
is bijective for any $D$. Using our lemma \ref{l3} one sees that it is enough to prove the following
lemma.

\begin{lem}\label{lt3'}
Let $C$ be a cofibrant and $\mathbb{U}$-small dg-category and
$M$ a $\mathbb{V}$-combinatorial $C(k)_{\mathbb{V}}$-model category
satisfying the same assumption as in lemma \ref{lt3}.
\begin{enumerate}
\item
Then, the Quillen adjunction
$$l_{!} : M^{C} \longrightarrow M^{\widehat{C}_{pe}} \qquad
M^{C} \longleftarrow M^{\widehat{C}_{pe}} : l^{*}$$
is a Quillen equivalence.
\item For any $F\in M^{\widehat{C}_{pe}}$, and
any a $\mathbb{U}$-small diagram of perfect and cofibrant objects in
$C^{op}-Mod_{\mathbb{U}}$,
$X : I \longrightarrow C^{op}-Mod_{\mathbb{U}}$,
the natural morphism
$$Hocolim_{i}F(X_{i}) \longrightarrow
F(Hocolim_{i}X_{i})$$
is an isomorphism in $Ho(M)$.
\item For any $F\in  M^{\widehat{C}_{pe}}$, and any
perfect complex $Z\in C(k)_{\mathbb{U}}$ and any $X\in M$, the natural
morphism
$$Z\otimes^{\mathbb{L}}F(X) \longrightarrow F(Z\otimes^{\mathbb{L}}X)$$
is an isomorphism in $Ho(M)$.
\end{enumerate}
\end{lem}

\textit{Proof:} This is the same as for lemma \ref{lt3}
and sub-lemma \ref{sllt3}. \hfill $\Box$ \\

Coming back to our square of dg-categories one sees that
the horizontal morphisms are quasi-fully faithful and that the right vertical morphism
is a quasi-equivalence. This formally implies that the left vertical morphism
is quasi-fully faithful. We now consider the square of sets
$$\xymatrix{
[\widehat{C}_{pe},\widehat{D}_{pe}] \ar[r] \ar[d] &
[\widehat{C}_{pe},\widehat{D}] \ar[d] \\
[C,\widehat{D}_{pe}] \ar[r] &  [C,\widehat{D}],}$$
obtained from the square of dg-categories by passing to equivalence classes
of objects. Again, the right vertical morphism is a bijection and the horizontal
morphisms are injective.
For $u\in [C,\widehat{D}_{pe}]$, its image in $[C,\widehat{D}]$
comes from an element $v\in [\widehat{C}_{pe},\widehat{D}]$.
For any $x\in C$,
$v(l(x)) \in \widehat{D}$ is a perfect $D^{op}$-module, and thus
so is $v(Z\otimes^{\mathbb{L}}l(x))\simeq Z\otimes^{\mathbb{L}}v(l(x))$
for any perfect complex $Z$ of $k$-modules.
As any perfect $C^{op}$-module is constructed as a retract
of a finite homotopy colimit of objects of the form
$Z\otimes^{\mathbb{L}}l(x)$, we deduce that
$v(X)$ is a perfect $D^{op}$-module for any
$X\in \widehat{C }_{pe}$. Therefore, Cor. \ref{clmono} implies that
$v$ comes in fact from an element in $[\widehat{C}_{pe},\widehat{D}_{pe}]$.
This shows that
$[\widehat{C}_{pe},\widehat{D}_{pe}] \longrightarrow [C,\widehat{D}_{pe}]$ is
surjective, and thus
that
$$\mathbb{R}\underline{Hom}(\widehat{C}_{pe},\widehat{D}_{pe}) \longrightarrow
\mathbb{R}\underline{Hom}(C,\widehat{D}_{pe})$$
is quasi-essentially surjective. This finishes the proof of the theorem. \hfill $\Box$ \\

The following corollary is the promised derived version of Morita theory.

\begin{cor}\label{ct3}
Let $C$ and $D$ be two $\mathbb{U}$-small dg-categories, then there exists a natural
isomorphism in $Ho(dg-Cat_{\mathbb{V}})$
$$\mathbb{R}\underline{Hom}_{c}(\widehat{C},\widehat{D})\simeq
\widehat{C^{op}\otimes^{\mathbb{L}}D}\simeq
Int((C\otimes^{\mathbb{L}}D^{op})-Mod_{\mathbb{U}}).$$
In particular, there exists a natural weak equivalence
$$Map_{c}(\widehat{C},\widehat{D})\simeq
|(C\otimes^{\mathbb{L}}D^{op})-Mod_{\mathbb{U}}|,$$
where $Map_{c}(\widehat{C},\widehat{D})$ is the sub-simplicial set
of continuous morphisms in $Map(\widehat{C},\widehat{D})$ and where
$|(C\otimes^{\mathbb{L}}D^{op})-Mod_{\mathbb{U}}|$ is the nerve
of the sub-category of equivalences between $C\otimes^{\mathbb{L}}D^{op}$-modules.
\end{cor}

\textit{Proof:} The first part follows from the universal properties
of internal Hom's, as by theorem \ref{t3}
$$\mathbb{R}\underline{Hom}_{c}(\widehat{C},\widehat{D})\simeq
\mathbb{R}\underline{Hom}(C,\mathbb{R}\underline{Hom}(D^{op},\widehat{\mathbf{1}}))\simeq
\mathbb{R}\underline{Hom}(C\otimes^{\mathbb{L}}D^{op},\widehat{\mathbf{1}})\simeq
\widehat{C^{op}\otimes^{\mathbb{L}}D}.$$
The second part follows from the relation between mapping spaces and
internal Hom's, as well as Prop. \ref{cp5'}. Indeed, one has
$$Map_{c}(\widehat{C},\widehat{D}) \simeq Map(\mathbf{1},\mathbb{R}\underline{Hom}_{c}(\widehat{C},\widehat{D}))\simeq
Map(\mathbf{1},\mathbb{R}\underline{Hom}(C\otimes^{\mathbb{L}}D^{op},\widehat{\mathbf{1}}))
\simeq
Map(\mathbf{1},\widehat{C^{op}\otimes^{\mathbb{L}}D}).$$
By theorem \ref{t1} this last simplicial set is weakly equivalent to
the nerve of the category of equivalences between
quasi-representable $\mathbb{V}$-small $\widehat{C^{op}\otimes^{\mathbb{L}}D}$-modules. The
enriched Yoneda
lemma for the model category $C\otimes^{\mathbb{L}}D^{op}-Mod$ easily implies that
this nerve is weakly equivalent to the nerve of equivalences
between $\mathbb{U}$-small $C\otimes^{\mathbb{L}}D^{op}$-modules. \hfill $\Box$ \\

\section{Applications}

In this last section we present three kinds of applications
of our main results. A first application explains the relation
between Hochschild cohomology and internal Hom's of dg-categories.
In the same spirit, we present a relation between the negative
part of Hochschild cohomology and the higher homotopy groups of
the \emph{classifying space of dg-categories}, as well as an
interpretation of the fundamental group of this space
as the so-called \emph{derived Picard group}.
As a second application, we present a proof of the existence of
a good localization functor for dg-categories. This implies for
example the existence of a quotient of a dg-category by
a full sub-dg-category, satisfying the required universal
property. Finally, our last application states that the
(derived) dg-category of morphisms between the dg-categories of quasi-coherent
complexes over some (reasonable) schemes is naturally equivalent to
the dg-category of quasi-coherent complexes over their product.
Under smoothness and properness conditions
the same statement stays correct
when one replaces \emph{quasi-coherent} by \emph{perfect}.
This last result can be considered as a solution to
a question of D. Orlov, concerning the existence of
representative objects for triangulated functors between
derived categories of smooth projective varieties.

\subsection{Hochschild cohomology, classifying space of dg-categories,
and derived Picard groups}

As a first application we give a formula relating higher homotopy groups
of mapping spaces between dg-categories and Hochschild cohomology. For this,
let us recall that for any $\mathbb{U}$-small dg-category $C$, one defines its
Hochschild cohomology groups as
$$\mathbb{HH}^{i}(C):=H^{i}(\mathbb{R}\underline{Hom}_{C\otimes^{\mathbb{L}}C^{op}}(C,C)),$$
where $C$ is the $C\otimes^{\mathbb{L}}C^{op}$-module defined by the trivial formula
$C(x,y):=C(x,y)$, and where $\mathbb{R}\underline{Hom}_{C\otimes^{\mathbb{L}}C^{op}}$ are the $Ho(C(k))$-enriched
Hom's of the category $Ho(C\otimes^{\mathbb{L}}C^{op}-Mod_{\mathbb{U}})$.
More generally, the Hochschild complex of $C$ is defined by
$$\mathbb{HH}(C):=\mathbb{R}\underline{Hom}_{C\otimes^{\mathbb{L}}C^{op}}(C,C)),$$
which is a well defined object in the derived category $Ho(C(k))$ of complexes
of $k$-modules.

\begin{cor}\label{chh-}
With the notation above, there exists an isomorphism in $Ho(C(k))$
$$\mathbb{HH}(C)\simeq \mathbb{R}\underline{Hom}(C,C)(Id,Id),$$
where $Id$ is the identity of $C$, considered as an
object of the dg-category $\mathbb{R}\underline{Hom}(C,C)$.
In particular, one has
$$\mathbb{HH}^{i}(C)\simeq H^{i}(\mathbb{R}\underline{Hom}(C,C)(Id,Id)).$$
\end{cor}

\textit{Proof:} Using Thm. \ref{t2}, one has
$$\mathbb{R}\underline{Hom}(C,C)(Id,Id)\simeq
Int(C\otimes^{\mathbb{L}}C^{op}-Mod_{\mathbb{U}}^{rqr}).$$
Furthermore, through this identification the identity morphism
of $C$ goes to the bi-module $C$ itself. This implies the result by
the definition of Hochschild cohomology.
\hfill $\Box$ \\

An important consequence of Cor. \ref{chh-} is the following Morita invariance
of Hochschild cohomology.

\begin{cor}\label{chh--}
With the notation above, there exists an isomorphism in $Ho(C(k))$
$$\mathbb{HH}(C)\simeq \mathbb{HH}(\widehat{C}).$$
\end{cor}

\textit{Proof:} Indeed, the identity of $\widehat{C}$ is clearly
continuous, and thus by Thm. \ref{t3} (1) one has
$$\mathbb{HH}(\widehat{C})\simeq
\mathbb{R}\underline{Hom}(\widehat{C},\widehat{C})(Id,Id)
\simeq \mathbb{R}\underline{Hom}(C,\widehat{C})(\underline{h},\underline{h}),$$
where $\underline{h} : C \longrightarrow \widehat{C}$ is the Yoneda
embedding. As the morphism $\underline{h}$ is quasi-fully faithful,
Cor. \ref{cp5'''} implies that the
morphism
$$\underline{h}^{*} : \mathbb{R}\underline{Hom}(C,\widehat{C})(\underline{h},\underline{h})
\longrightarrow \mathbb{R}\underline{Hom}(C,C)(Id,Id)$$
is a quasi-isomorphism. Cor. \ref{chh-} implies the result. \hfill $\Box$ \\

\begin{cor}\label{chh}
With the notation above one has isomorphisms of groups
$$\pi_{i}(Map(C,C),Id)\simeq \mathbb{HH}^{1-i}(C)$$
for any $i>1$. For $i=1$, one has an isomorphism of groups
$$\pi_{1}(Map(C,C),Id)\simeq \mathbb{HH}^{0}(C)^{*}=Aut_{Ho(C\otimes^{\mathbb{L}}C^{op}-Mod_{\mathbb{U}})}(C).$$
\end{cor}

\textit{Proof:} This follows immediately from Thm. \ref{t1},
the well-known relations between mapping spaces and classifying
spaces of model categories (see e.g. \cite[Cor. A.0.4]{hagII}) and the formula
$$H^{-i}(\mathbb{R}\underline{Hom}_{C\otimes^{\mathbb{L}}C^{op}}(C,C))\simeq
\pi_{i}(Map_{C\otimes^{\mathbb{L}}C^{op}-Mod_{\mathbb{U}}}(C,C)).$$
\hfill $\Box$ \\

Let $|dg-Cat|$ be the nerve of the category of quasi-equivalences in
$dg-Cat_{\mathbb{U}}$. Using the usual relations between
mapping spaces in model category and nerve of categories of equivalences
(see e.g. \cite[Appendix A]{hagII}) one finds the following consequence.

\begin{cor}\label{chh'}
For a $\mathbb{U}$-small dg-category $C$, one has
$$\pi_{i}(|dg-Cat|,C)\simeq \mathbb{HH}^{2-i}(C) \qquad \forall \; i >2.$$
Moreover, one has
$$\pi_{2}(|dg-Cat|,C)\simeq \mathbb{HH}^{0}(C)^{*}.$$
\end{cor}

\begin{rmk}
\emph{The above corollary only gives an interpretation
of negative Hochschild cohomology groups. The positive
part of the Hochschild cohomology can also be
interpreted in terms of deformation theory of
dg-categories as done for example in \cite[\S 8.5]{hagII}.}
\end{rmk}

For a ($\mathbb{U}$-small) dg-algebra $A$, one can define
the derived Picard group $RPic(A)$ of $A$, as done
for example in \cite{rz,ke2,yek}.
Using our notations and definitions, the group
$RPic(A)$ can be defined in the following way.
To simplify notations let us assume that
the underlying complex of $A$ is cofibrant, and we will
consider $A$ as a dg-category with a unique object
which we denote by $BA$.
Note that the category
$(A\otimes A^{op})-Mod_{\mathbb{U}}$,
of $A\otimes A^{op}$-dg-modules, is also
the category $(BA\otimes BA^{op})-Mod_{\mathbb{U}}$.
This category can be
endowed with the following
monoidal structure. For $X$ and $Y$ two
$(A\otimes A^{op})$-dg-modules,
we can form the internal tensor product
$X\otimes_{A}Y \in (A\otimes A^{op})-Mod_{\mathbb{U}}$
as the coequalizer of the two natural morphisms
$$(X\otimes A \otimes Y) \rightrightarrows X\otimes Y.$$
This endows the model category
$(A\otimes A^{op})-Mod_{\mathbb{U}}$ with
a structure of monoidal model category (see for
example \cite{kt} where the simplicial analog
is considered). Deriving this monoidal
structure provides a monoidal category
$(Ho((A\otimes A^{op})-Mod_{\mathbb{U}}),\otimes_{A}^{\mathbb{L}})$.
By definition,
the group $RPic(A)$ is the group of isomorphism
classes of objects in $Ho((A\otimes A^{op})-Mod_{\mathbb{U}})$
which are invertible for the monoidal structure $\otimes_{A}^{\mathbb{L}}$.

\begin{cor}\label{crpic}
There is
a group isomorphism
$$RPic(A)\simeq \pi_{1}(|dg-Cat_{\mathbb{V}}|,\widehat{BA}).$$
\end{cor}

\textit{Proof:} This easily follows from the formula
$$\pi_{1}(|dg-Cat_{\mathbb{V}}|,\widehat{C}) \simeq
Aut_{Ho(dg-Cat)}(\widehat{C})$$
and Cor. \ref{ct3}. \hfill $\Box$ \\

\subsection{Localization and quotient of dg-categories}

Let $C$ be a $\mathbb{U}$-small dg category, and $S$ be a set of
morphisms in $[C]$. For any $\mathbb{U}$-small dg-category $D$, we consider
$Map_{S}(C,D)$ the sub-simplicial set of
$Map(C,D)$ being the union of all connected components corresponding to
morphisms $f : C \longrightarrow D$ in $Ho(dg-Cat)$ such that
$[f] : [C] \longrightarrow [D]$ sends $S$ to isomorphisms in $[D]$.

\begin{cor}\label{cloc}
The $Ho(SSet_{\mathbb{U}})$-enriched functor
$$Map_{S}(C,-) : Ho(dg-Cat_{\mathbb{U}}) \longrightarrow
Ho(SSet_{\mathbb{U}})$$
is co-represented by an object $L_{S}(C) \in Ho(dg-Cat_{\mathbb{U}})$.
\end{cor}

\textit{Proof:} Let $I_{k}$ be the dg-category with two objects $0$ and $1$, and
freely generated by a unique morphism $0 \rightarrow 1$. Using
theorem \ref{t1} one easily sees that $Map(I_{k},C)$ can be identified
with the nerve of the category $(C^{op}-Mod_{\mathbb{U}})_{rqr}^{I}$, of morphisms
between quasi-representable $C^{op}$-modules. Using the dg-Yoneda
lemma one sees that $[I_{k},C]$ is in a natural bijection with
isomorphism classes of morphisms in $[C]$. In particular, the set
$S$ can be classified by a morphism in $Ho(dg-Cat_{\mathbb{U}})$
$$S : \coprod_{f\in S}I_{k} \longrightarrow C.$$
We consider the natural morphism $I_{k} \longrightarrow \mathbf{1}$,
and we define $L_{S}C$ to be the homotopy push-out
$$\xymatrix{
\coprod_{f\in S}I_{k} \ar[r] \ar[d] & C \ar[d] \\
\coprod_{f\in S}\mathbf{1} \ar[r] & L_{S}C.}$$

For any $D$ one has a homotopy pull-back diagram
$$\xymatrix{
Map(L_{S}C,D) \ar[r] \ar[d] & \prod_{f\in S}Map(\mathbf{1},D) \ar[d] \\
Map(C,D) \ar[r] & \prod_{f\in S}Map(I_{k},D).}$$
Therefore, in order to see that $L_{S}C$ has the correct universal property, it is
enough to check that $Map(\mathbf{1},D) \longrightarrow Map(I_{k},D)$ induces an injection
on $\pi_{0}$, a bijection on $\pi_{i}$ for $i>0$, and that its image in
$[I_{k},D]$ consists of all morphisms in $[D]$ which are isomorphisms.
Using theorem \ref{t1} once again we see that this follows from the following
very general fact: if $M$ is a model category, then
the Quillen adjunction $Mor(M) \leftrightharpoons M$ (where $Mor(M)$ is the model
category of morphisms in $M$), sending a morphism in $M$ to its
target, induces a fully faithful functor $Ho(M) \longrightarrow Ho(Mor(M))$, whose essential image
consists of all equivalences in $M$. \hfill $\Box$.

\begin{cor}\label{cloc2}
Let $C\in dg-Cat_{\mathbb{U}}$ be a
dg-category and $S$ a set of morphisms in $[C]$.
Then, the natural morphism $C \longrightarrow L_{S}C$
induces for any $D\in dg-Cat_{\mathbb{U}}$
a quasi-fully faithful morphism
$$\mathbb{R}\underline{Hom}(L_{S}C,D) \longrightarrow
\mathbb{R}\underline{Hom}(C,D),$$
whose quasi-essential image consists
of all morphisms $C \rightarrow D$ in $Ho(dg-Cat)$
sending $S$ to isomorphisms in $[D]$.
\end{cor}

\textit{Proof:} This follows formally from Cor. \ref{cloc},
Thm. \ref{t2} and Lem. \ref{lmono}. \hfill $\Box$ \\

One important example of application of the localization construction
is the existence of a good theory of quotients of dg-categories.
For this, let $C$ be a $\mathbb{U}$-small dg-category, and
$\{X_{i}\}_{i\in I}$ be a sub-set of objects in $C$.
We assume that $[C]$ has a zero object $0$.
One consider $S$ the set of morphisms in $[C]$ consisting
of all $X_{i} \rightarrow 0$. The dg-category
$L_{S}C$ is then denoted by $C/<X_{i}>$, and is called
the quotient of $C$ by the sub-set of objects $\{X_{i}\}_{i\in I}$.
This terminology is justified by the fact that for any
dg-category $D$ with a zero object, the morphism
$$l^{*} : \mathbb{R}\underline{Hom}(C/<X_{i}>,D) \longrightarrow
\mathbb{R}\underline{Hom}(C,D)$$
is quasi-fully faithful, and its image consists of all
morphisms $f : C \longrightarrow D$ such that for all $i\in I$
$[f(X_{i})]\simeq 0$ in $[D]$.

\subsection{Maps between dg-categories of quasi-coherent complexes}

We now pass to our last application describing maps between
dg-categories of quasi-coherent complexes on $k$-schemes.
For this, let $X$
be a quasi-compact and separated scheme over $Spec\, k$. We consider
$QCoh(X)$  the category of $\mathbb{U}$-small quasi-coherent
sheaves on $X$. As this is a Grothendieck category we know that
there exists a $\mathbb{U}$-cofibrantly generated model category
$C(QCoh(X))$  of (unbounded) complexes of quasi-coherent sheaves
on $X$ (the cofibrations being the monomorphisms and the
equivalences being the quasi-isomorphisms, see e.g. \cite{ho2}).
It is easy to check that the natural
$C(k)_{\mathbb{U}}$-enrichment of $C(QCoh(X))$ makes it into a
$C(k)_{\mathbb{U}}$-model category, and thus as explained in
\S 3 we can construct a $\mathbb{V}$-small dg-category
$Int(C(QCoh(X))$. This dg-category will be denoted by
$L_{qcoh}(X)$. Note that $[L_{qcoh}(X)]$ is naturally equivalent to 
the (unbounded) 
derived category of quasi-coherent sheaves $D_{qcoh}(X)$, and will be
identified with it.

We need to recall that an object $E$ in $L_{qcoh}(X)$
is homotopically finitely presented, or perfect in the sense of \S 7,
if and only if it is a compact object of $D_{qcoh}(X)$, and thus
if and only if it is a perfect complex on $X$ (see for example
\cite{bv}). We will use this fact implicitly in the sequel.
\\

\begin{thm}\label{tfour}
Let $X$ and $Y$ be two quasi-compact and separated schemes over $k$,
and assume that one of them is flat over $Spec\, k$.
Then, there
exists an isomorphism in $Ho(dg-Cat_{\mathbb{V}})$
$$\mathbb{R}\underline{Hom}_{c}(L_{qcoh}(X),L_{qcoh}(Y))\simeq
L_{qcoh}(X\times_{k} Y).$$
\end{thm}

\textit{Proof:} We start noticing that the model
categories $C(QCoh(X))$ and $C(QCoh(Y))$ are stable,
proper, cofibrantly generated,
and admit a compact generator (see \cite{bv}). Therefore, they satisfy the conditions of the
main theorem of \cite{ss}, and thus one can find two objects
$E_{X}$ and $E_{Y}$ in $L_{qcoh}(X)$ and $L_{qcoh}(Y)$, and two Quillen
equivalences
$$C(QCoh(X)) \leftrightharpoons A_{X}^{op}-Mod_{\mathbb{U}} \qquad
C(QCoh(Y)) \leftrightharpoons A_{Y}^{op}-Mod_{\mathbb{U}}$$
where
$A_{X}$ (resp. $A_{Y}$) is the full sub-dg-category of
$L_{qcoh}(X)$ (resp. of $L_{qcoh}(Y)$) consisting of $E_{X}$
(resp. $E_{Y}$) only (in other words, $A_{X}$ is the dg-category with a unique
object and $\mathbb{R}\underline{End}(E_{X})$ as endomorphism dg-algebra).
In the following we will write $A_{X}$ for both, the dg-category
and the corresponding dg-algebra $\mathbb{R}\underline{End}(E_{X})$
(and the same with $A_{Y}$).
These Quillen equivalences
are $C(k)$-enriched Quillen equivalences, and with
a bit of care one can check that they
provide natural isomorphisms in $Ho(dg-Cat_{\mathbb{V}})$
$$L_{qcoh}(X)\simeq \widehat{A_{X}} \qquad L_{qcoh}(Y)\simeq \widehat{A_{Y}}.$$

\begin{lem}\label{ldual}
There exists an isomorphism in $Ho(dg-Cat_{\mathbb{V}})$
$$\widehat{A_{Y}}\simeq \widehat{A_{Y}^{op}}.$$
\end{lem}

\textit{Proof:} By the general theory of \cite{ss} it is enough to show that
the triangulated category $D_{qcoh}(Y)\simeq [L_{qcoh}(Y)]$ possesses
a compact generator $F_{Y}$ such that
the dg-algebra $\mathbb{R}\underline{End}(F_{Y})$ is
naturally equivalent to $\mathbb{R}\underline{End}(E_{Y})^{op}$.
For this we take $F_{Y}=E_{Y}^{\vee}$ to be the dual perfect complex
of $E_{Y}$.
Let $<F_{Y}>$ be the smallest thick triangulated sub-category
of $D_{parf}(Y)$ containing $F_{Y}$. We let
$\phi : D_{parf}(Y) \longrightarrow D_{parf}(Y)^{op}$ be the
involution sending a perfect complex $E$ to its dual $E^{\vee}$.
Then, clearly $\phi(<F_{Y}>)=<E_{Y}>=D_{parf}(Y)$.
This shows that $F_{Y}$ classically generates
$D_{parf}(Y)$, and thus by \cite[Thm. 2.1.2]{bv} that
$F_{Y}$ is a compact generator of $D_{qcoh}(Y)$.  \hfill $\Box$ \\

\begin{lem}\label{ldiag}
There exists an isomorphism in $Ho(dg-Cat_{\mathbb{V}})$
$$\widehat{A_{X}\otimes_{k}^{\mathbb{L}} A_{Y}}\simeq L_{qcoh}(X\times_{k} Y).$$
\end{lem}

\textit{Proof:} This follows from the fact that the external product
$E_{X}\boxtimes E_{Y}$ is a compact generator of
$D_{qcoh}(X\times_{k} Y)$, as explained in \cite[Lem. 3.4.1]{bv}. Indeed,
flat base change induces a natural quasi-isomorphism of dg-algebras
(one uses here that either $X$ or $Y$ is flat over $k$)
$$\mathbb{R}\underline{End}(E_{X}\boxtimes E_{Y})\simeq
\mathbb{R}\underline{End}(E_{X})\otimes_{k}^{\mathbb{L}} \mathbb{R}\underline{End}(E_{Y})\simeq
A_{X}\otimes_{k}^{\mathbb{L}} A_{Y}.$$
\hfill $\Box$ \\

We are now ready to prove theorem \ref{tfour}. Indeed, using theorem \ref{t3} one finds
$$\mathbb{R}\underline{Hom}_{c}(L_{qcoh}(X),L_{qcoh}(Y))\simeq
\mathbb{R}\underline{Hom}_{c}(\widehat{A_{X}},\widehat{A_{Y}})\simeq
\mathbb{R}\underline{Hom}(A_{X},\widehat{A_{Y}}).$$
Lemma \ref{ldual} and the universal properties of internal Hom's give an isomorphism
$$\mathbb{R}\underline{Hom}(A_{X},\widehat{A_{Y}})\simeq
\mathbb{R}\underline{Hom}(A_{X},\widehat{A_{Y}^{op}})\simeq
\widehat{A_{X}\otimes_{k}^{\mathbb{L}} A_{Y}}.$$
Finally lemma \ref{ldiag} implies the theorem. \hfill $\Box$ \\

\begin{cor}\label{cfour-1}
Under the same conditions as in Thm. \ref{tfour},
there
exists a bijection
between $[L_{qcoh}(X),L_{qcoh}(Y)]_{c}$, the sub-set of
$[L_{qcoh}(X),L_{qcoh}(Y)]$ consisting of continuous morphisms,
and the isomorphism classes of objects in the
derived category $D_{qcoh}(X\times_{k} Y)$.
\end{cor}

\textit{Proof:} Readily follows from theorem \ref{tfour} and the fact that
$[L_{qcoh}(X\times_{k} Y)]\simeq D_{qcoh}(X\times_{k} Y)$. \hfill $\Box$ \\

Tracking back the construction of the equivalence in theorem
\ref{tfour} one sees that the bijection of corollary 
\ref{cfour-1} can be described as follows. Let 
$E \in D_{qcoh}(X\times_{k} Y)$ be an object, and let 
us consider the two projections
$$p_{X} : X\times_{k}Y \longrightarrow X \qquad 
p_{Y} : X\times_{k}Y \longrightarrow Y.$$
We consider the functor
$$\phi_{E} : D_{qcoh}(X) \longrightarrow D_{qcoh}(Y)$$
defined by 
$$\phi_{E}(F):=\mathbb{R}(p_{Y})_{*}(\mathbb{L}p_{X}^{*}(F)\otimes^{\mathbb{L}}
E),$$
for any $F\in D_{qcoh}(X)$. Then, the functor $\phi_{E}$ is the natural 
functor induced by the morphism $L_{qcoh}(X) \longrightarrow L_{qcoh}(Y)$
in $Ho(dg-Cat)$, corresponding to $E$ via the bijection 
of Cor. \ref{cfour-1}. 

\begin{cor}\label{cfour}
Let $X$ be a quasi-compact and separated scheme, flat over $Spec\, k$.
Then, one has
$$\pi_{1}(Map(L_{qcoh}(X),L_{qcoh}(X),Id))\simeq \mathcal{O}_{X}(X)^{*}$$
$$\pi_{i}(Map(L_{qcoh}(X),L_{qcoh}(X),Id))\simeq \mathbb{HH}^{1-i}(A_{X}) \simeq 0 \; \forall \; i>1.$$
\end{cor}

\textit{Proof:} Indeed theorem \ref{tfour}, theorem \ref{t1}, corollary
\ref{cint} and corollary \ref{cp5'} give
$$Map(L_{qcoh}(X),L_{qcoh}(X))\simeq Map(*,L_{qcoh}(X\times_{k} X)).$$
Furthermore, the identity on the right is clearly sent to the
diagonal $\Delta_{X}$ in $L_{qcoh}(X\times_{k} X)$.
Therefore, one finds for any $i>1$
$$\pi_{i}(Map(L_{qcoh}(X),L_{qcoh}(X)),Id))\simeq
\pi_{i}(Map(*,L_{qcoh}(X\times_{k} X)),\Delta_{X})\simeq$$
$$H^{1-i}(L_{qcoh}(X\times_{k} X)(\Delta_{X},\Delta_{X}))\simeq
Ext^{1-i}_{X\times_{k} X}(\Delta(X),\Delta(X))\simeq
0.$$
For $i=1$, one has
$$\pi_{1}(Map(L_{qcoh}(X),L_{qcoh}(X)),Id))\simeq
\pi_{1}(Map(*,L_{qcoh}(X\times_{k} X)),\Delta_{X})\simeq
 Aut_{D_{qcoh}(X\times_{k} X)}(\Delta_{X})\simeq \mathcal{O}_{X}(X)^{*}.$$
\hfill $\Box$ \\

Corollary \ref{cfour}
combined with the usual relations between mapping spaces and nerves
of categories of equivalences also has the following important consequence.

\begin{cor}\label{cfour'}
Let $X$ be a quasi-compact and separated scheme, flat over $k$. Then, one has
$$\pi_{i}(|dg-Cat|,L_{qcoh}(X))\simeq 0 \qquad \forall \; i>2.$$
In particular, the sub-simplicial set of $|dg-Cat|$ corresponding to dg-categories
of the form $L_{qcoh}(X)$, for $X$ a quasi-compact and separated scheme
flat over $k$, is a $2$-truncated simplicial set.
\end{cor}

We finish by a refined version of theorem \ref{tfour} involving
only perfect complexes instead of all quasi-coherent complexes.
For this, we will denote by $L_{parf}(X)$ the full sub-dg-category of
$L_{qcoh}(X)$ consisting of all perfect complexes.

\begin{thm}\label{tfour2}
Let $X$ and $Y$ be two smooth and proper schemes over $Spec\, k$.
Then, there
exists an isomorphism in $Ho(dg-Cat_{\mathbb{V}})$
$$\mathbb{R}\underline{Hom}(L_{parf}(X),L_{parf}(Y))\simeq
L_{parf}(X\times_{k} Y).$$
\end{thm}

\textit{Proof:} The triangulated category
$D_{qcoh}(X)$ being generated by its compact objects, one sees
that the Yoneda embedding
$$L_{qcoh}(X) \longrightarrow \widehat{L_{parf}(X)}$$
is an isomorphism in $Ho(dg-Cat_{\mathbb{V}})$. Using
our Thm. \ref{t3} we see that
$\mathbb{R}\underline{Hom}(L_{parf}(X),L_{parf}(Y))$ can be
identified, up to quasi-equivalence, with the full sub-dg-category
of $\mathbb{R}\underline{Hom}_{c}(L_{qcoh}(X),L_{qcoh}(Y))$
consisting of all morphisms $L_{qcoh}(X) \longrightarrow
L_{qcoh}(Y)$ which preserve perfect complexes. Using Thm. \ref{tfour},
we see that $\mathbb{R}\underline{Hom}(L_{parf}(X),L_{parf}(Y))$
is quasi-equivalent to the full sub-dg-category
of $L_{qcoh}(X\times_{k}Y)$ consisting of objects
$E$ such that for any perfect complex $F$ on $X$, the
complex $\mathbb{R}(p_{Y})_{*}(p_{X}^{*}(F)\otimes^{\mathbb{L}}E)$
is perfect on $Y$. To finish the proof we thus need to
show that an object $E\in D_{qcoh}(X\times_{k}Y)$ is perfect
if and only if the functor
$$\Phi_{E}:=\mathbb{R}(p_{Y})_{*}(p_{X}^{*}(-)\otimes^{\mathbb{L}}E) :
D_{qcoh}(X) \longrightarrow D_{qcoh}(Y)$$
preserves perfect objects. Clearly, as $X$ is flat and proper
over $Spec\, k$, $\Phi_{E}$ preserves perfect complexes if
$E$ is itself perfect.

Conversely, let $E$ be an object in $D_{qcoh}(X\times_{k}Y)$
such that $\Phi_{E}$ preserves perfect complexes.

\begin{lem}\label{ltfour2}
Let $Z$ be a smooth and proper scheme over $Spec\, k$, and
$E\in D_{qcoh}(Z)$. If for any perfect complex $F$ on $Z$, the
complex of $k$-modules $\mathbb{R}\underline{Hom}(F,E)$ is
perfect, then $E$ is perfect on $Z$.
\end{lem}

\textit{Proof of the lemma:} We let $A_{Z}$ be a dg-algebra
over $k$ such that $L_{qcoh}(Z)$ is quasi-equivalent to
$\widehat{A_{Z}}$ (with the same abuse of notations that
$A_{Z}$ also means the dg-category with a unique object
and $A_{Z}$ as endomorphism dg-algebra). As $Z$ is flat and
proper over $Spec\, k$, the underlying complex of $k$-modules
of $A_{Z}$ is perfect. Furthermore, as $Z$ is smooth, the diagonal
$\Delta : Z \hookrightarrow Z\times_{k}Z$ is a local complete intersection
morphism, and thus $\Delta_{*}(\mathcal{O}_{Z})$ is a perfect complex
on $Z$. Equivalently, the $A_{Z}\otimes^{\mathbb{L}}_{k}A_{Z}^{op}$-dg-module
$A_{Z}$ is perfect, or equivalently lies in the smallest
sub-dg-category of $\widehat{A_{Z}^{op}\otimes^{\mathbb{L}}_{k}A_{Z}}$
containing $A_{Z}\otimes^{\mathbb{L}}A_{Z}^{op}$ and which is
stable by retracts, homotopy push-outs and the loop functor
(or the shift functor).

We now apply our theorem \ref{t3} in order to
translate this last fact in terms of dg-categories of morphisms.
Let $F : \widehat{A_{Z}} \longrightarrow \widehat{A_{Z}}$
be the morphism of dg-categories sending
an $A_{Z}^{op}$-dg-module $M$ to the free $A_{Z}^{op}$-dg-module
$$F(M):=\underline{M}\otimes^{\mathbb{L}}A_{Z}^{op},$$
where $\underline{M}$ is the underlying complex of
$k$-modules of $M$. By what we have seen, the identity morphism lies in the
smallest sub-dg-category of $\mathbb{R}\underline{Hom}_{c}(\widehat{A_{Z}},\widehat{A_{Z}})$ containing
the object $F$ and which is stable by retracts, homotopy push-outs
and the loop functor.
Evaluating the identity of the dg-category
$\widehat{A_{Z}}$ at an object $M$, we get that the
object $M\in \widehat{A_{Z}}$ lies in the smallest sub-dg-category
of $\widehat{A_{Z}}$ containing $\underline{M}\otimes^{\mathbb{L}} A_{Z}^{op}$
and stable by retracts, homotopy push-outs and by the loop functor. 
Now, by our hypothesis
the object $E$ corresponds to $M\in \widehat{A_{Z}}$ such that
$\underline{M}$ is a perfect complex
of $k$-modules. Therefore, $M$ itself belongs to the smallest
sub-dg-category of $\widehat{A_{Z}}$ containing $A_{Z}^{op}$ and
which stable by retracts, homotopy push-outs and the loop functor. 
By definition
of being perfect, this implies that $M\in \widehat{A_{Z}}_{pe}$, and thus
that $E$ is a perfect complex on $Z$. \hfill $\Box$ \\

Let now $E_{X}$ and $E_{Y}$ be compact generators of $D_{qcoh}(X)$
and $D_{qcoh}(Y)$. Then, by the projection formula one has
$$\mathbb{R}(p_{Y})_{*}((E_{X}^{\vee}\boxtimes E_{Y}^{\vee})\otimes^{\mathbb{L}}E)
\simeq \mathbb{R}(p_{Y})_{*}(p_{X}^{*}(E_{X}^{\vee})\otimes^{\mathbb{L}}E)
\otimes^{\mathbb{L}}E_{Y}^{\vee}$$
which is perfect on $Y$. This implies in particular that 
$$\mathbb{R}\underline{Hom}(E_{X}\boxtimes E_{Y},E)\simeq
\mathbb{R}\Gamma(Y,\mathbb{R}(p_{Y})_{*}((E_{X}^{\vee}\boxtimes E_{Y}^{\vee})\otimes^{\mathbb{L}}E)$$
is a perfect complex of $k$-modules. 
As the perfect complex $E_{X}\boxtimes E_{Y}$
is a generator on $D_{qcoh}(X\times_{k}Y)$, 
one sees that
for any perfect complex $F$ on $X\times_{k}Y$, the complex
of $k$-modules $\mathbb{R}\underline{Hom}(F,E)$ is perfect.
The lemma \ref{ltfour2} implies that $E$ is perfect on 
$X\times_{k}Y$, which
finishes the proof of the theorem. \hfill $\Box$ \\

\end{document}